\documentclass[11pt, reqno]{amsart}
\usepackage{amsthm,amsmath,amssymb}
\usepackage{mathtools}
\usepackage{graphicx}
\usepackage{caption}
\usepackage{ragged2e}
\graphicspath{{Figure/}}
\numberwithin{equation}{section}
\numberwithin{figure}{section}
\theoremstyle{plain}

\newtheorem{thm}{Theorem}[section]
\newtheorem{lem}[thm]{Lemma}

\theoremstyle{definition}
\theoremstyle{remark}

\begin{document}
\setlength{\abovedisplayskip}{10pt}
\setlength{\belowdisplayskip}{10pt}

\title{Lozenge Tilings of Hexagons with Intrusions I: Generalized Intrusion}

\author{Seok Hyun Byun}
\address{School of Mathematical and Statistical Sciences, Clemson University, Clemson, SC 29634, U.S.A.}
\email{sbyun@clemson.edu}
\thanks{S.B. is supported by the AMS-Simons Travel Grant.}

\author{Tri Lai}
\address{Department of Mathematics, University of Nebraska – Lincoln, Lincoln, NE 68588, U.S.A.}
\email{tlai3@unl.edu}
\thanks{T.L. was supported in part by Simons Collaboration Grant (\#585923).}

\begin{abstract}
MacMahon's classical theorem on the number of boxed plane partitions has been generalized in several directions. One way to generalize the theorem is to view boxed plane partitions as lozenge tilings of a hexagonal region and then generalize it by making some holes in the region and counting its tilings. In this paper, we provide new regions whose numbers of lozenges tilings are given by simple product formulas. The regions we consider can be obtained from hexagons by removing structures called \textit{intrusions}. In fact, we show that the tiling generating functions of those regions under certain weights are given by similar formulas. These give the $q$-analogue of the enumeration results.
\end{abstract}

\maketitle

\section{Introduction}

The enumeration of lozenge tilings has received considerable attention during the last three decades. Most of the works in this area were motivated by the classical result of MacMahon~\cite{M}. Via the bijection of David and Tomei~\cite{DT}, MacMahon's theorem on boxed plane partitions can be rephrased as follows. The number of lozenge tilings of hexagons with sides of length $x,y,z,x,y,$ and $z$ (in clockwise order) is given by the following product formula:

\begin{equation}\label{eqn:eq11}
    \prod_{i=1}^{x}\prod_{j=1}^{y}\prod_{k=1}^{z}\frac{i+j+k-1}{i+j+k-2}=\frac{H(x)H(y)H(z)H(x+y+z)}{H(x+y)H(y+z)H(z+x)},
\end{equation}
where $H(n):=\prod_{i=0}^{n-1}i!$ for $n\geq 1$ and $H(0):=1$.

\begin{figure}
    \centering
    \includegraphics[width=0.6\textwidth]{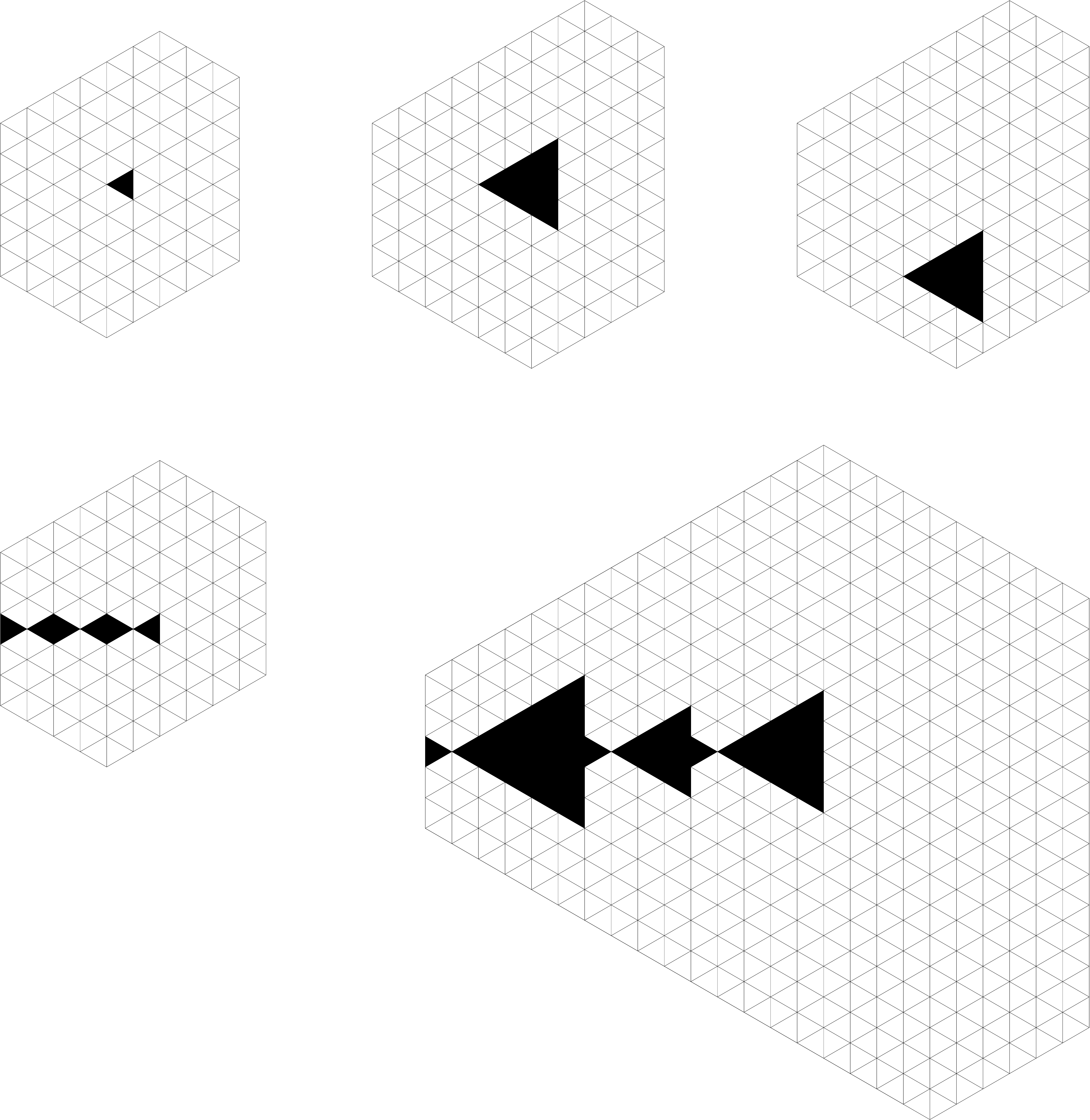}
    \caption{Various regions obtained from hexagons by making some holes in them.}
    \label{fig:Figure_1}
\end{figure}

From a lozenge tiling point of view, one way to generalize MacMahon's theorem is to make some holes in the hexagonal region and to count its tilings (see ~\cite{B,C2, C3, CEKZ, CK3, CL, Fu2, L1, L2, L3, L4, LR1, R} for various results that generalized MacMahon's theorem in this direction.). In 1998, Krattenthaler and Okada~\cite{KO} enumerated the lozenge tilings of a hexagon with the central unit triangle removed. Soon after, in 2001, this result was generalized by Ciucu et al.~\cite{CEKZ}. They generalized Krattenthaler and Okada's result by considering a hexagonal region with a triangle of arbitrary size removed from the center. This more general result was then further generalized by Rosengren. In his paper~\cite{R}, Rosengren considered a hexagon with a triangular hole in an arbitrary position and found a formula that enumerated its lozenge tilings. In fact, he proved a more general statement. He considered certain weights on lozenges and found a tiling generating function of the region under the weight assignment (see the three pictures at the top of Figure \ref{fig:Figure_1} that illustrate these regions).

Recently, the first author~\cite{B} considered a hexagonal region with several unit triangles removed\footnote{In that paper, the collection of these removed unit triangles is called ``an intrusion". In this paper, the definition of ``an intrusion" is more general than that.} and enumerated the lozenge tilings of the region (see the bottom left picture in Figure \ref{fig:Figure_1} for an example). Motivated by the result of Ciucu et al.~\cite{CEKZ} and Rosengren~\cite{R}, we generalize this result by allowing removed left-pointing triangles to have arbitrary odd size (except the rightmost one, which also allows having an arbitrary even size) and considering tiling generating functions under a similar weight assignment (see the bottom right picture in Figure \ref{fig:Figure_1} that illustrates these regions). A detailed description of the regions is given in the next section. Recently, Fulmek~\cite{Fu2} also considered lozenge tilings enumeration of hexagonal regions that generalize the region considered in \cite{B}. While Fulmek considered new regions by translating the intrusions considered in ~\cite{B}, we generalize the intrusions in ~\cite{B} by inflating some left-pointing unit triangles. We also obtained the generalization of the work of the first author in different directions, and they will be presented in the forthcoming papers~\cite{BL1,BL2}.

This paper is organized as follows. 
In Section 2, we first state a weighted generalization of MacMahon's theorem. Then, we describe our new regions and present the main result.
In Section 3, we gather several results that are needed in the next section.
In Section 4, we give a proof of the main theorem, Theorem 2.2.

\section{New regions and statement of the main results}

\begin{figure}
    \centering
    \includegraphics[width=0.55\textwidth]{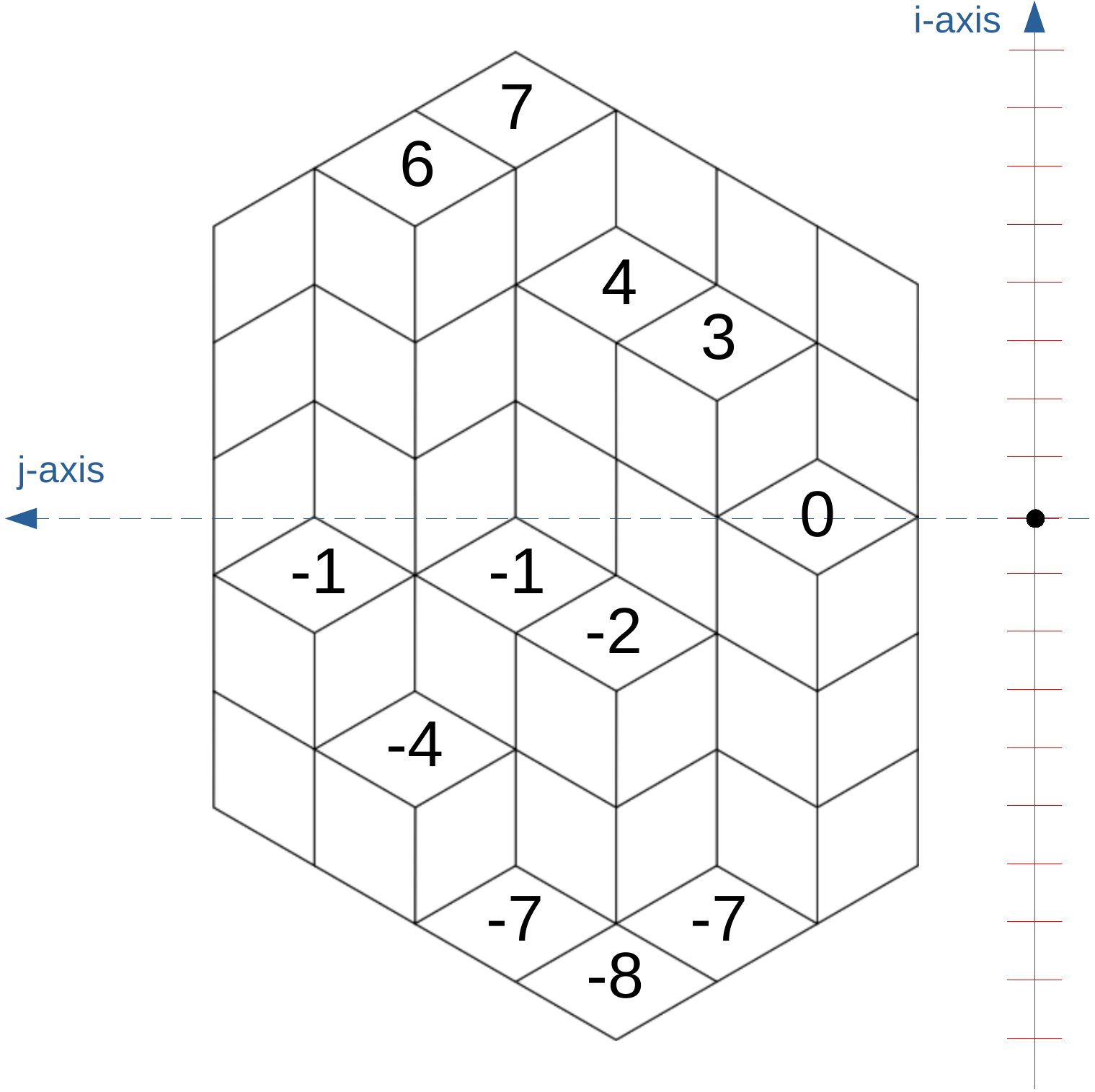}
    \caption{A lozenge tiling of the region $H_{5,3,4}$ and the weight assignment on its lozenges. Horizontal lozenges labeled by an integer $n$ are given a weight of $\frac{q^{n}+q^{-n}}{2}$ and all other lozenges are given a weight of 1.}
    \label{fig:Figure_2}
\end{figure}

We consider a triangular lattice such that one family of lattice lines is vertical. On the lattice, we draw a hexagon with sides of length $x,y,z,x,y,$ and $z$ clockwise from the left and denote the region by $H_{x,y,z}$. We first provide a weighted generalization of MacMahon's theorem since it is needed to understand the main result better. To give weights on lozenges, we put the region $H_{x,y,z}$ on $(i,j)-$coordinate system\footnote{We set the unit length of the $i$-axis on this plane equals the half of the side length of unit triangles. We use this $(i,j)-$coordinate system throughout this paper.} as indicated in Figure \ref{fig:Figure_2}. As shown in the picture, we put the region so that the $j-$axis overlaps with the perpendicular bisector of the left side. Now, we give a weight $\frac{q^i+q^{-i}}{2}$ to all horizontal lozenges whose centers have $i-$coordinate $i$ and give a weight $1$ to all nonhorizontal lozenges. \textit{A lozenge tiling} of a region is a collection of unit lozenges that covers the region without gaps or overlaps. When certain weights are assigned on lozenges and a tiling of a region is given, \textit{a weight of the tiling} of the region is the product of the weights of all the lozenges that constitute the tiling. \textit{A tiling generating function (TGF)} of the region is the sum of weights of all its lozenge tilings. Under this weight assignment, let $M_{q}(H_{x,y,z})$ be the TGF of $H_{x,y,z}$. To state the lemma, we introduce several notations. Let $\langle n\rangle_q:=\frac{q^n-q^{-n}}{q-q^{-1}}$ for positive integers $n$. Since $\langle n\rangle_{q}\rightarrow n$ as $q\rightarrow 1$, $\langle n\rangle_{q}$ can be considered as a $q$-analogue of a positive integer $n$. Also, we set $\langle n\rangle_q!:=\prod_{i=1}^{n}\langle i\rangle_{q}$ for positive integers $n$ and $\langle0\rangle_q!:=1$. The following lemma is not new. This is equivalent to a special case of Rosengren's result (Theorem 2.1 in ~\cite{R}) when a triangle of size $0$ is removed from a hexagonal region. The second author also proved a result (Lemma 2.7 in ~\cite{L5}) that can be specialized to almost the same result as the following lemma\footnote{In \cite{L5}, the second author found a TGF of dented semihexagons (which can be specialized to hexagons by making suitable choices of positions of dents) with a different choice of $j$-axis. However, the proof presented in that paper can also be employed to give the proof of Lemma 2.1 in the current paper. A more general result that generalizes both Lemma 2.1 in this paper and Lemma 2.7 in ~\cite{L5}, together with its proof will be presented in the forthcoming paper of the authors \cite{BL1}.}.

\begin{lem}
For nonnegative integers $x,y,$ and $z,$ the TGF of the region $H_{x,y,z}$ is given by the following formula: \begin{equation}\label{eqn:eq21}
        M_q(H_{x,y,z})=\Bigg[\prod_{i=1}^{y}\prod_{j=1}^{z}\frac{q^{i-j}+q^{j-i}}{2}\Bigg]\cdot\frac{H_{q}(x)H_{q}(y)H_{q}(z)H_{q}(x+y+x)}{H_{q}(x+y)H_{q}(y+z)H_{q}(z+x)},
    \end{equation}
where $H_{q}(n):=\prod_{i=0}^{n-1}\langle i\rangle_{q}!$ for positive integers $n$ and $H_{q}(0):=1$.
\end{lem}

We generalize this region $H_{x,y,z}$ by putting some triangular holes in it. For any nonnegative integers $a,m,n,k$, and $a_1,\ldots, a_k$, we set\footnote{When $k=0$, we set $\{a_1,\ldots,a_k\}=\emptyset$ and thus $\{b_1,\ldots,b_k\}=\emptyset$.} $b_i=\sum_{j=i}^{k}a_j$ for $i\in\{1,\ldots,k\}$. We now define two regions $H_{e}(a,m,n,k:a_1,\ldots,a_k)$ and $H_{o}(a,m,n,k:a_1,\ldots,a_k)$.

\begin{figure}
    \centering
    \includegraphics[width=1.0\textwidth]{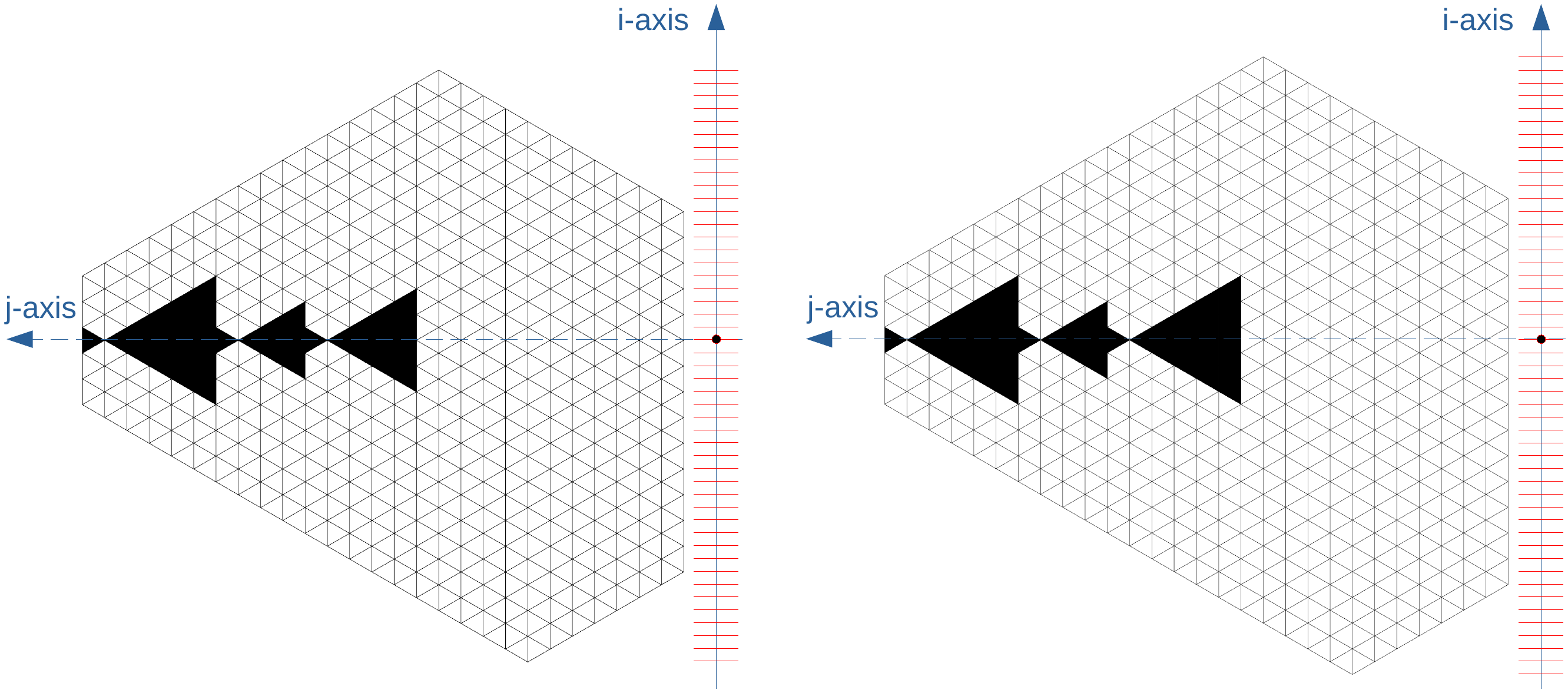}
    \caption{Two regions $H_{e}(2,4,8,3:2,1,2)$ (left) and $H_{o}(2,4,8,3:2,1,2)$ (right) on the $(i,j)-$coordinate planes.}
    \label{fig:Figure_3}
\end{figure}

When $k=0$, we set $H_{e}(a,m,n,0:\cdot):=H_{2a,m,n}$. When $k\geq 1$, we first consider a hexagon with sides of length $2a+1, m+2b_1+k-1, n+k, 2a+2b_1, m+k,$ and $n+2b_1+k-1$ clockwise from the left. Then, we delete the consecutive triangles of side lengths $1, 2a_1+1, 1, 2a_2+1,\ldots,1, 2a_{k-1}+1, 1, 2a_{k}$ on the perpendicular bisector of the left side from the left so that they alternate the orientation, where the first triangle is pointing right. We denote this region by $H_{e}(a,m,n,k:a_1,\ldots,a_k)$ (see the left picture in Figure \ref{fig:Figure_3}). The collection of $2k$ triangles we removed is an \textit{intrusion} of the region.

Similarly, when $k=0$, we set $H_{o}(a,m,n,0:\cdot):=H_{2a+1, m, n}$. When $k\geq1$, we consider a hexagon with sides of length $2a+1, m+2b_1+k, n+k, 2a+2b_1+1, m+k$, and $n+2b_1+k$ clockwise from the left. Then, we delete the consecutive triangles of side lengths $1, 2a_1+1, 1, 2a_2+1,\ldots,1, 2a_{k-1}+1, 1, 2a_k+1$ on the perpendicular bisector of the left side from the left so that they alternate the orientation, where the first one is pointing right. We denote this region by $H_{o}(a,m,n,k:a_1,\ldots,a_k)$ (see the right picture in Figure \ref{fig:Figure_3}). Again, an \textit{intrusion} of this region is the collection of $2k$ triangles that we removed.

We also assign certain weights to lozenges. Following the weight assignment on lozenges in $H_{x,y,z}$, we placed the regions $H_{e}(a,m,n,k:a_1,\ldots,a_k)$ and $H_{o}(a,m,n,k:a_1,\ldots,a_k)$ on the $(i,j)-$plane as indicated in Figure \ref{fig:Figure_3}. As before, we put the regions so that the $j-$axis overlaps with the perpendicular bisectors of the left sides. Then we give a weight $\frac{q^i+q^{-i}}{2}$ to all horizontal lozenges whose centers have $i-$coordinate $i$ and give a weight $1$ to all nonhorizontal lozenges. Under this weight assignment, let $M_q(H_{e}(a,m,n,k:a_1,\ldots,a_k))$ and $M_q(H_{o}(a,m,n,k:a_1,\ldots,a_k))$ be the TGFs of the two regions. Note that if we interchange $m$ and $n$, the TGFs are invariant. This is because we assign weights to the lozenges so that it is symmetric with respect to the $j-$axis. Therefore, without loss of generality, we assume $n\geq m$ in the theorem below.

\begin{thm}\label{thm:Theorem1}
For any nonnegative integers $a,m,n,k$, and $a_1,\ldots,a_k$ such that $n\geq m$,

\begin{equation}\label{eqn:eq22}
\begin{aligned}
    &\frac{M_q(H_{e}(a,m,n,k:a_1,\ldots,a_k))}{M_q(H_{2a+2b_1, m+k, n+k})}\\
    &
\begin{aligned}
    =\prod_{i=0}^{k-1}\Bigg[&\frac{(\langle b_{1}+i+\frac{1}{2}\rangle_{q^2})_{m+k-b_{1}-2i}(\langle a+b_1-b_{i+1}+i+1\rangle_{q^2})_{m+k+2b_{i+1}-2i-1}}{(\langle a+b_1+i+\frac{1}{2}\rangle_{q^2})_{m+k-2i}(\langle i+1\rangle_{q^2})_{m+k+b_{i+1}-2i-1}}\\
    &\cdot \frac{(\langle a+n+k+b_1-i\rangle_{q^2})_{b_{i+1}}}{2^{4b_{i+1}}(\langle i+1\rangle_{q^2})_{b_1}(\langle n+k-i-\frac{1}{2}\rangle_{q^2})_{b_{i+1}}(\langle a+b_1+i+1\rangle_{q^2})_{-b_{i+1}}}\\
    &\cdot \frac{(\langle m+k-i+\frac{1}{2}\rangle_{q^2})_{\lfloor\frac{n-m}{2}\rfloor}(\langle n+k-i\rangle_{q^2})_{-\lfloor\frac{n-m}{2}\rfloor}}{(\langle a+m+k+b_{1}-i+\frac{1}{2}\rangle_{q^2})_{\lfloor\frac{n-m}{2}\rfloor}(\langle a+n+k+b_{1}-i\rangle_{q^2})_{-\lfloor\frac{n-m}{2}\rfloor}}\\
    &\cdot \frac{(\langle m+k-i\rangle_{q^2})_{\lfloor\frac{n-m}{2}\rfloor}(\langle n+k-i-\frac{1}{2}\rangle_{q^2})_{-\lfloor\frac{n-m}{2}\rfloor}}{(\langle m+k+b_{i+1}-i\rangle_{q^2})_{\lfloor\frac{n-m}{2}\rfloor}(\langle n+k+b_{i+1}-i-\frac{1}{2}\rangle_{q^2})_{-\lfloor\frac{n-m}{2}\rfloor}}\\
    &\cdot \prod_{j=1}^{i}\Big[
    \frac{(\langle b_j-b_{i+1}+i+1-j\rangle_{q^2})^2}{(\langle i+1-j\rangle_{q^2})^2}\cdot(\langle b_{j+1}+i-j+\frac{1}{2}\rangle_{q^2})_{a_j}(\langle b_{j+1}+i-j+1\rangle_{q^2})_{a_j}\Big]\Bigg]
\end{aligned}
\end{aligned}
\end{equation}
and
\begin{equation}\label{eqn:eq23}
\begin{aligned}
    &\frac{M_q(H_{o}(a,m,n,k:a_1,\ldots,a_k))}{M_q(H_{2a+2b_1+1, m+k, n+k})}\\
    &
\begin{aligned}
    =\prod_{i=0}^{k-1}\Bigg[&\frac{(\langle a+b_1-b_{i+1}+i+1\rangle_{q^2})_{n+k+2b_{i+1}-2i}(\langle b_1+i+\frac{3}{2}\rangle_{q^2})_{n+k-b_1-2i-2}}{(\langle i+1\rangle_{q^2})_{n+k+b_{i+1}-2i-1}(\langle a+b_1+i+\frac{3}{2}\rangle_{q^2})_{n+k-2i-1}}\\
    &\cdot \frac{(\langle b_{i+1}+1\rangle_{q^2})(\langle a+m+k+b_1-i+1\rangle_{q^2})_{b_{i+1}}}{2^{4b_{i+1}+2}(\langle i+2\rangle_{q^2})_{b_1}(\langle m+k-i+\frac{1}{2}\rangle_{q^2})_{b_{i+1}}(\langle a+b_1+i+1\rangle_{q^2})_{-b_{i+1}}}\\
    &\cdot \frac{(\langle m+k-i\rangle_{q^2})_{\lfloor\frac{n-m}{2}\rfloor}(\langle n+k-i-\frac{1}{2}\rangle_{q^2})_{-\lfloor\frac{n-m}{2}\rfloor}}{(\langle a+m+k+b_1-i+1\rangle_{q^2})_{\lfloor\frac{n-m}{2}\rfloor}(\langle a+n+k+b_1-i+\frac{1}{2}\rangle_{q^2})_{-\lfloor\frac{n-m}{2}\rfloor}}\\
    &\cdot \frac{(\langle m+k-i+\frac{1}{2}\rangle_{q^2})_{\lfloor\frac{n-m}{2}\rfloor}(\langle n+k-i\rangle_{q^2})_{-\lfloor\frac{n-m}{2}\rfloor}}{(\langle m+k+b_{i+1}-i+\frac{1}{2}\rangle_{q^2})_{\lfloor\frac{n-m}{2}\rfloor}(\langle n+k+b_{i+1}-i\rangle_{q^2})_{-\lfloor\frac{n-m}{2}\rfloor}}\\
    &
    \begin{aligned}
        \cdot \prod_{j=1}^{i}&\Big[\frac{(\langle b_j-b_{i+1}+i+1-j\rangle_{q^2})^2}{(\langle i+1-j\rangle_{q^2})(\langle b_j+i+1-j\rangle_{q^2})}\frac{(\langle b_j+i+2-j\rangle_{q^2})}{(\langle i+2-j\rangle_{q^2})}\\
        &\cdot (\langle b_{j+1}+i-j+\frac{3}{2}\rangle_{q^2})_{a_j}(\langle b_{j+1}+i-j+2\rangle_{q^2})_{a_j}\Big]\Bigg],
    \end{aligned}
\end{aligned}
\end{aligned}
\end{equation}
where $\langle n\rangle_q:=\frac{q^n-q^{-n}}{q-q^{-1}}$ for positive integers (or half-integers) $n$ and 

\begin{equation*}
    (\langle a\rangle_{q})_n:=
    \begin{cases}
        \prod_{i=0}^{n-1}\langle a+i\rangle_{q}  & \text{if $n$ is positive,}\\
        1                    & \text{if $n$ is $0$,}\\
        \frac{1}{\prod_{i=1}^{-n}\langle a-i\rangle_{q}} & \text{if $n$ is negative,}
    \end{cases}
\end{equation*}
for positive integers (or half-integers) $a$ and integers $n$ such that $a+n>0$.
\end{thm}

Note that by combining \eqref{eqn:eq21} and \eqref{eqn:eq22} (or \eqref{eqn:eq23}), one can obtain a product formula for the tiling generating function of the region $H_{e}(a,m,n,k:a_1,\ldots,a_k)$ (or $H_{o}(a,m,n,k:a_1,\ldots,a_k)$). Hence, one can view our results as tiling generating functions of these two regions. The case when $q=0$ and $k=0$ gives MacMahon's formula~\cite{M} and the case when $q=0$ and $a_{i}=0$ for all $i$ was considered by the first author~\cite{B}.

\section{Some previous results}

In this section, we gather some previous results that we utilize in the next section where we give a proof of Theorem 2.2. We first recall the result of Lai and Rohatgi \cite{LR2} that gives TGFs of two families of regions $R_{\mathbf{l},n,x}$ and $\overline{R}_{\mathbf{l},n,x}$ that we define below. In their notation, $\mathbf{l}=\{l_1,\ldots,l_m\}$ is a set of positive integers or an empty set\footnote{When $\mathbf{l}=\emptyset$, we set $m=0$ and $l_{0}=0$ by convention.}, where elements are written in increasing order, $n$ is a nonnegative integer, and $x$ is an integer that makes all the side lengths of the regions nonnegative.

\begin{figure}
    \centering
    \includegraphics[width=16.5cm]{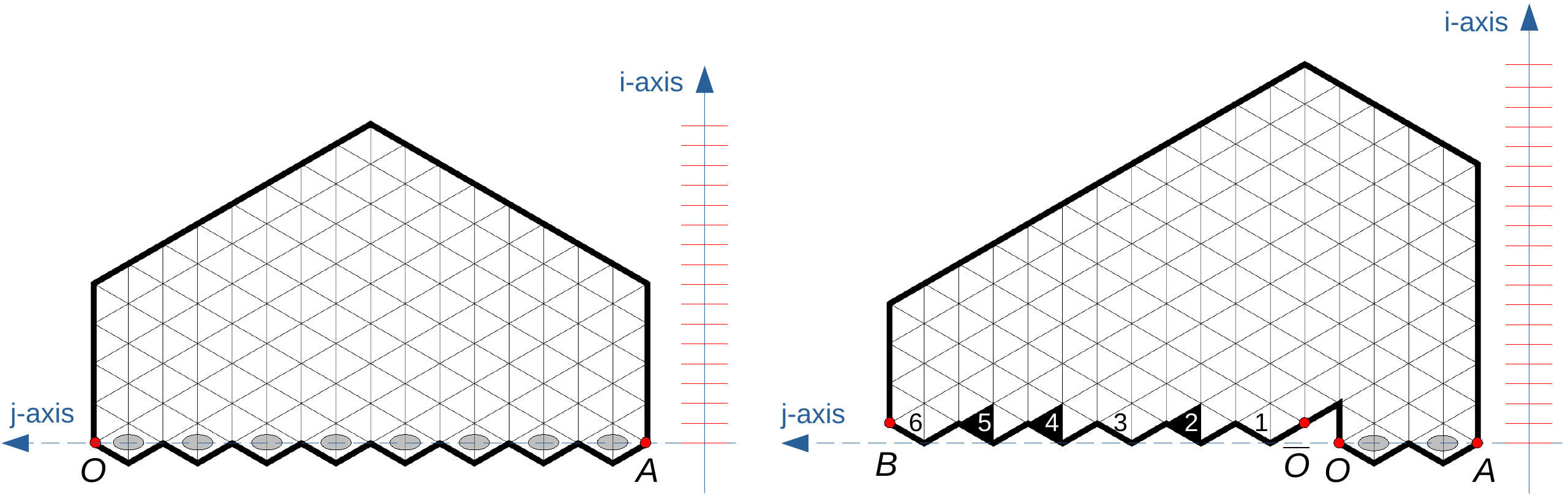}
    \caption{The regions $R_{\emptyset,8,3}$ (left) and $R_{\{1,3,6\},2,3}$ (right) on the $(i,j)-$planes.}
\end{figure}
\begin{figure}
    \centering
    \includegraphics[width=16.5cm]{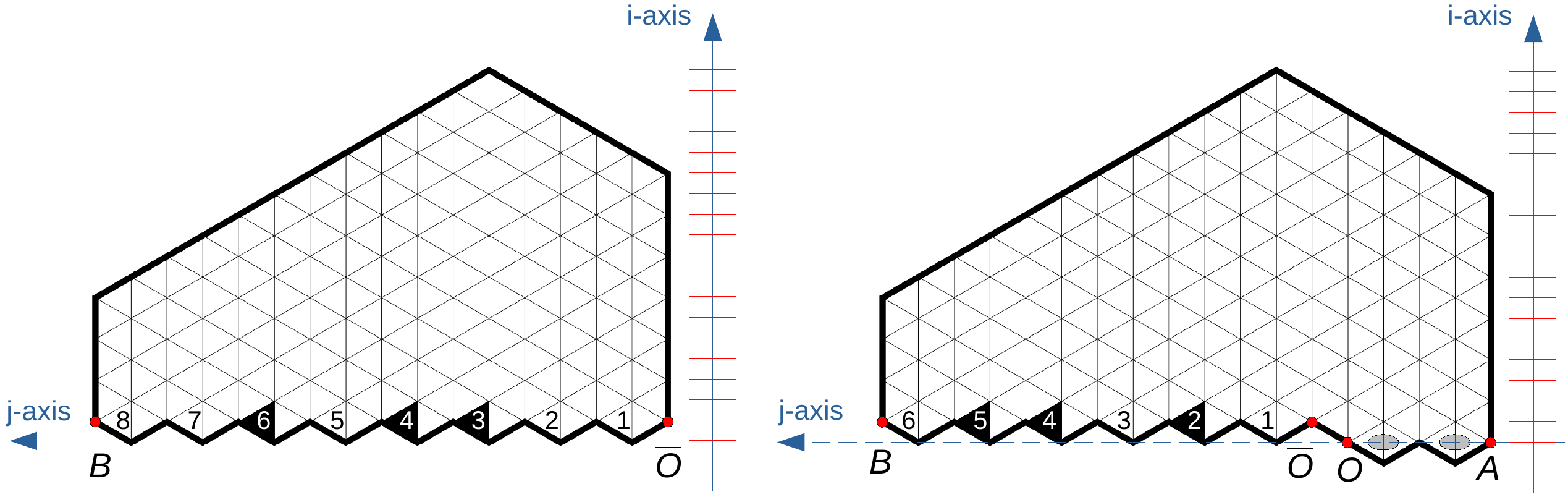}
    \caption{The regions $\overline{R}_{\mathbf\{1,2,5,7,8\},0,3}$ (left) and $\overline{R}_{\mathbf\{1,3,6\},2,3}$ (right) on the $(i,j)-$planes.}
\end{figure}

First, we define the region $R_{\mathbf{l},n,x}$. Let \emph{O} be any lattice point on a triangular lattice, and \emph{$\Bar{O}$} be a lattice point that is one unit length away from \emph{O} in the northwest direction. Consider a horizontal line \emph{l} through \emph{O} (which is not a lattice line) and let \emph{A} be the $n$th lattice point to the right of \emph{O} which is on \emph{l} (when $n=0$, we set $A=O$). Similarly, consider the horizontal line $\Bar{\emph{l}}$ through \emph{$\Bar{O}$} and let \emph{B} be the $l_m$th lattice point to the left of \emph{$\Bar{O}$} which is on $\Bar{\emph{l}}$ (when $\mathbf{l}=\emptyset$, we set $B=\Bar{O}$). Then, we label left-pointing unit triangles on $\Bar{\emph{l}}$ to the left of \emph{$\Bar{O}$} by $1,2,3,\ldots$ from the right to the left. The region $R_{\mathbf{l},n,x}$ is defined differently, depending on whether $\mathbf{l}$ is an empty set or not (see Figure 3.1 for two examples).

When $\mathbf{l}=\emptyset$, the region $R_{\emptyset,n,x}$ is defined as follows without using \emph{$\Bar{O}$} and \emph{B}. From \emph{A} to \emph{O}, we follow the zigzag line along \emph{l}, by alternating moving one unit to the southwest and one unit to the northwest. Then, from \emph{O}, we move $x+1$ units to the north, $n$ units to the northeast, $n$ units to the southeast, and then $x+1$ units to the south until we reach \emph{A}. For any integer $x\geq-1$, $R_{\emptyset,n,x}$ is defined to be the bounded region enclosed by the path described above.

When $\mathbf{l}\neq\emptyset$, from \emph{A} to \emph{O}, we follow the zigzag line along \emph{l}, as we did in the previous case. Then, from \emph{O} to \emph{$\Bar{O}$}, we move one unit to the north and then one unit to the southwest. Now, we connect \emph{$\Bar{O}$} and \emph{B} using the same type of zigzag line along $\Bar{\emph{l}}$. From \emph{B}, we move $x$ units to the north, $2l_m-m+n+1$ units to the northeast, $m+n$ units to the southeast, and then $x+l_m-m+1$ units to the south until we reach \emph{A}. From the bounded region enclosed by the path described above, we delete the labeled left-pointing unit triangles on $\Bar{\emph{l}}$ whose labels are in $[l_m]\setminus\mathbf{l}$.\footnote{For any positive integer $k$, $[k]:=\{1,2,\ldots,k\}$.} For any nonnegative integer $x$, $R_{\mathbf{l},n,x}$ is defined to be the resulting region.

The region $\overline{R}_{\mathbf{l},n,x}$ is defined differently, depending on whether $n$ is $0$ or not (see Figure 3.2). Its construction is almost identical to $R_{\mathbf{l},n,x}$. The differences are as follows.

1) We connect \emph{O} and \emph{$\Bar{O}$} using the unit segment that connects them.

2-1) When $n=0$, from \emph{B}, we move $x$ units to the north, $2l_m-m$ units to the northeast, $m$ units to the southeast, and then $x+l_m-m$ units to the south until we reach \emph{$\Bar{O}$}. Then, we removed the labeled left-pointing unit triangles on $\Bar{\emph{l}}$ whose labels are in $[l_m]\setminus\mathbf{l}$. For any nonnegative integer $x$, we denote the resulting region by $\overline{R}_{\mathbf{l},0,x}$ (so \emph{O} and \emph{A} are not used when we define the region $\overline{R}_{\mathbf{l},0,x}$).

2-2) When $n\geq1$, from \emph{B}, we move $x$ units to the north, $2l_m-m+n$ units to the northeast, $m+n+1$ units to the southeast, and then $x+l_m-m$ units to the south until we reach \emph{A}. Then, we remove the labeled left-pointing unit triangles on $\Bar{\emph{l}}$ whose labels are in $[l_m]\setminus\mathbf{l}$. For any nonnegative integer $x$, we denote the resulting regions by $\overline{R}_{\mathbf{l},n,x}$.

To assign weights on lozenges, we put the regions $R_{\mathbf{l},n,x}$ and $\overline{R}_{\mathbf{l},n,x}$ on the $(i,j)-$coordinate system as indicated in Figures 3.1 and 3.2. Then we assign weights on lozenges using the same criteria as before (giving weight $\frac{q^i+q^{-i}}{2}$ to all horizontal lozenges whose centers have $i-$coordinate $i$ and giving weight $1$ to all nonhorizontal lozenges) except the lozenges on the $j$-axis. According to the criteria, those lozenges on the $j$-axis are given a weight of $\frac{q^{0}+q^{0}}{2}=1$. However, instead of giving a weight $1$, we assign a weight $\frac{1}{2}$ to these lozenges. The lozenges weighted by $\frac{1}{2}$ are marked with shaded ellipses on the pictures in this paper. Under this weight assignment, let $M_{q}(R_{\mathbf{l},n,x})$ and $M_{q}(\overline{R}_{\mathbf{l},n,x})$ be the TGFs of the regions $R_{\mathbf{l},n,x}$ and $\overline{R}_{\mathbf{l},n,x}$, respectively. Ciucu first studied the lozenge tilings of these regions in ~\cite{
C2}. Recently, Lai and Rohatgi ~\cite{LR2} found the TGFs of these regions under this weight assignment\footnote{The regions $R_{\mathbf{l},n,x}$ and $\overline{R}_{\mathbf{l},n,x}$ correspond to the regions $T_{x,n,l_{m}}(\mathbf{l};[n])$ and $S_{x,n,l_{m}}(\mathbf{l};[n])$ in their paper.}.

\begin{thm} [\cite{LR2}, Theorems 2.12 and 2.13]
    The TGFs of the regions $R_{\mathbf{l},n,x}$ and $\overline{R}_{\mathbf{l},n,x}$ are given by the following product formulas:\footnote{The formulas presented in this paper looks different from the original formulas of Lai and Rohatgi and it is because two papers are using different $q$-analogue of integers. While they used $[n]_{q}\coloneqq\frac{1-q^{n}}{1-q}$, we are using $\langle n\rangle_{q}\coloneqq\frac{q^{n}-q^{-n}}{q-q^{-1}}$.}
    \begin{equation}
    \begin{aligned}
        M_{q}(R_{\mathbf{l},n,x})=&2^{-e}\frac{\prod_{1\leq i<j\leq m}\langle2(l_{j}-l_{i})\rangle_{q}\prod_{1\leq i<j\leq n}\langle2(j-i)\rangle_{q}}{\prod_{i=1}^{m}\prod_{j=1}^{n}\langle2(l_{i}+j)\rangle_{q}\prod_{i=1}^{m}\langle2l_{i}\rangle_{q}!\prod_{j=1}^{n}\langle2j-1\rangle_{q}!}\\
        &\cdot\prod_{i=1}^{\lceil n/2\rceil}\prod_{j=1}^{2n-4i+3}\langle2x+2l_{m}-2m+2i+j\rangle_{q}\\
        &\cdot\prod_{i=1}^{m}\langle2(x+l_{m}-m+n+i)\rangle_{q}\langle2(x+l_{m}-m+n+i+1)\rangle_{q}\\
        &\cdot\prod_{i=1}^{m}\prod_{j=1}^{n+i-1}\langle2x+2l_{m}-2m+n+i+j+1\rangle_{q}\\
        &\cdot\prod_{i=1}^{n}\prod_{j=1}^{m}\frac{\langle2(x+l_{m}-m+i+j-1)\rangle_{q}}{\langle2(x+l_{m}-m+i+j-1)+1\rangle_{q}}\\
        &\cdot\prod_{i=1}^{m}\prod_{j=1}^{l_{i}-i}\langle2(x+l_{m}-m+n+i+j+1)\rangle_{q}\langle2(x+l_{m}-i-j+1)\rangle_{q}
    \end{aligned}
    \end{equation}
    and
    \begin{equation}
    \begin{aligned}
        M_{q}(\overline{R}_{\mathbf{l},n,x})=&2^{-\overline{e}}\frac{\prod_{1\leq i<j\leq m}\langle2(l_{j}-l_{i})\rangle_{q}\prod_{1\leq i<j\leq n}\langle2(j-i)\rangle_{q}}{\prod_{i=1}^{m}\prod_{j=1}^{n}\langle2(l_{i}+j)\rangle_{q}\prod_{i=1}^{m}\langle2l_{i}-1\rangle_{q}!\prod_{j=1}^{n}\langle2j\rangle_{q}!}\\
        &\cdot\prod_{i=1}^{\lceil m/2\rceil}\prod_{j=1}^{2m-4i+3}\langle2x+2l_{m}-2m+2i+j-1\rangle_{q}\\
        &\cdot\prod_{i=1}^{n}\langle2(x+l_{m}+i)\rangle_{q}\prod_{i=1}^{n}\prod_{j=1}^{m+i}\langle2x+2l_{m}-m+i+j\rangle_{q}\\
        &\cdot\prod_{i=1}^{m}\prod_{j=1}^{n}\frac{\langle2(x+l_{m}-m+i+j-1)\rangle_{q}}{\langle2(x+l_{m}-m+i+j-1)+1\rangle_{q}}\\
        &\cdot\prod_{i=1}^{m}\prod_{j=1}^{l_{i}-i}\langle2(x+l_{m}-m+i+j+n)\rangle_{q}\langle2(x+l_{m}-i-j+1)\rangle_{q}
    \end{aligned}
    \end{equation}
where
$e=\frac{n(n+1)}{2}+(n+1)m+\sum^{m}_{i=1}(2l_i-i)$ \text{and} $\overline{e}=\frac{n(n+1)}{2}+n(m+1)+\sum^{m}_{i=1}(2l_i-i)$.\\
Also, when $\mathbf{l}=\emptyset$, we set $m=0$ and $l_{0}=0$ by convention.
\end{thm}

We also recall the Matching Factorization Theorem of Ciucu~\cite{C1}. Instead of stating the theorem in full generality, we present how the theorem is applied to the regions $H_{e}(a,m,m,k:a_1,\ldots,a_k)$ and $H_{o}(a,m,m,k:a_1,\ldots,a_k)$. Note that the region $H_{e}(a,m,m,k:a_1,\ldots,a_k)$ is symmetric and has a horizontal symmetry axis, which is the $j$-axis. Furthermore, since the symmetry axis is the $j$-axis and $\frac{q^{i}+q^{-i}}{2}=\frac{q^{-i}+q^{i}}{2}$, the weight assignment on the lozenges in this region is symmetric (in other words, two lozenges that are in symmetric positions have the same weight). On this region, we consider the midpoint of the vertical side of the deleted triangle of side length $2a_k$ (the rightmost triangle of the intrusion). Starting from the midpoint, we follow the zigzag line along the symmetry axis, alternating its direction between the northeast and the southeast (see the left picture in Figure 3.3). The region is then divided into two subregions by the zigzag line. Note that all lozenges in these two regions have weights that were inherited from the weight assignment of lozenges in $H_{e}(a,m,m,k:a_1,\ldots,a_k)$. We denote the top region by $H_{e}^{-}(a,m,m,k:a_1,\ldots,a_k)$ and let $M_{q}(H_{e}^{-}(a,m,m,k:a_1,\ldots,a_k))$ be the TGF of the region under the inherited weights on its lozenges. The region at the bottom is denoted by $H_{e}^{+}(a,m,m,k:a_1,\ldots,a_k)$. We also consider its TGF, but we first modify the weights on some of its lozenges. Note that the horizontal lozenges adjacent to the zigzag line were originally weighted by $\frac{q^{0}+q^{0}}{2}=1$, since they were on the $j$-axis. We replace the weights on these horizontal lozenges by $\frac{1}{2}$ and keep the weights of the other lozenges the same. Under this modified weight assignment, let $M_{q}(H_{e}^{+}(a,m,m,k:a_1,\ldots,a_k))$ be the TGF of $H_{e}^{+}(a,m,m,k:a_1,\ldots,a_k)$. Then, the Matching Factorization Theorem provides the following factorization of the TGF of $H_{e}(a,m,m,k:a_1,\ldots,a_k)$:

\begin{figure}
    \centering
    \includegraphics[width=16.5cm]{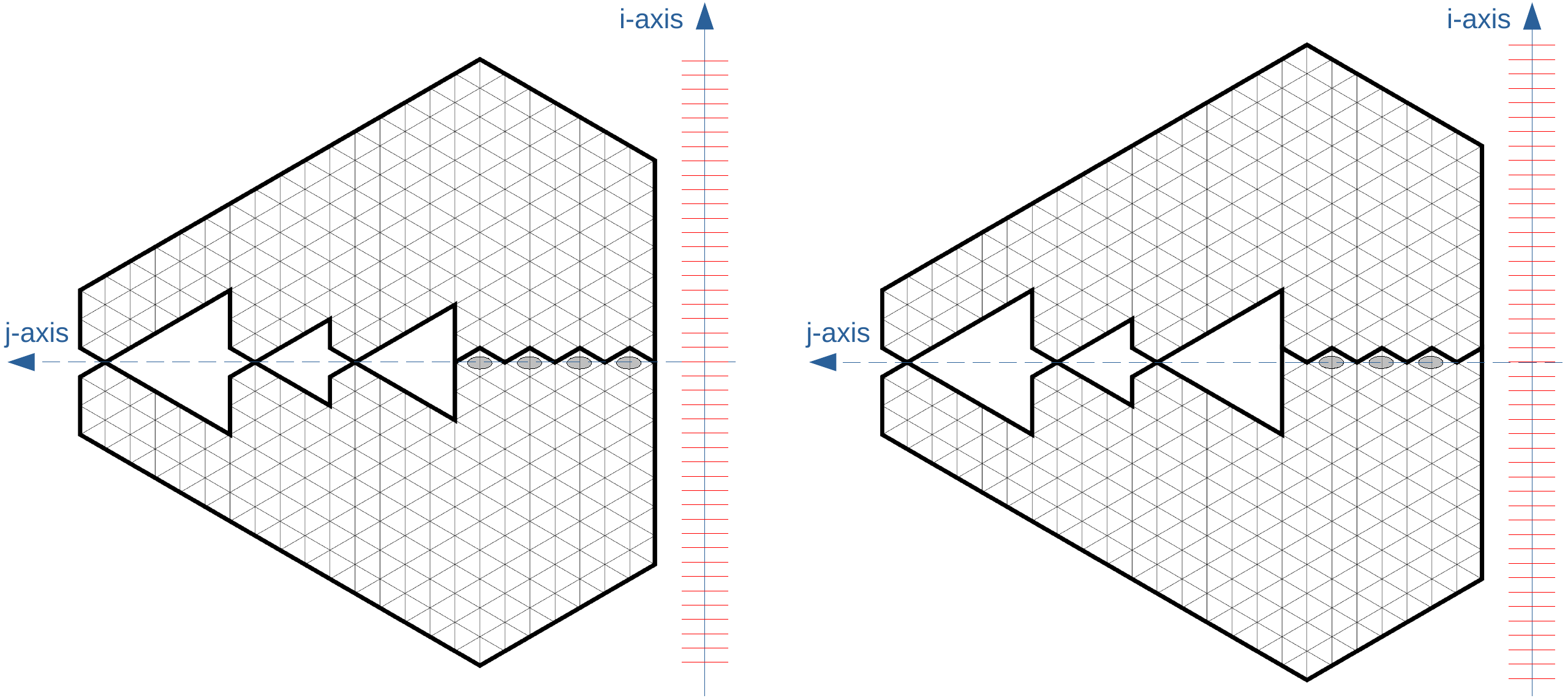}
    \caption{Pictures that show how the Factorization Theorem is applied on the regions $H_{e}(2,4,4,3:2,1,2)$ and $H_{o}(2,4,4,3:2,1,2)$. Left picture shows how the Factorization Theorem split the region $H_{e}(2,4,4,3:2,1,2)$ into two subregions $H_{e}^{+}(2,4,4,3:2,1,2)$ (bottom) and $H_{e}^{-}(2,4,4,3:2,1,2)$ (top). The picture on the right presents $H_{o}^{+}(2,4,4,3:2,1,2)$ (bottom) and $H_{o}^{-}(2,4,4,3:2,1,2)$ (top) that can be obtained from $H_{o}(2,4,4,3:2,1,2)$ by applying the theorem.}
\end{figure}

\begin{equation}
\begin{aligned}
    &M_{q}(H_{e}(a,m,m,k:a_1,\ldots,a_k))\\
    &=2^{m}M_{q}(H_{e}^{+}(a,m,m,k:a_1,\ldots,a_k))M_{q}(H_{e}^{-}(a,m,m,k:a_1,\ldots,a_k)).
\end{aligned}
\end{equation}

\begin{figure}
    \centering
    \includegraphics[width=16.5cm]{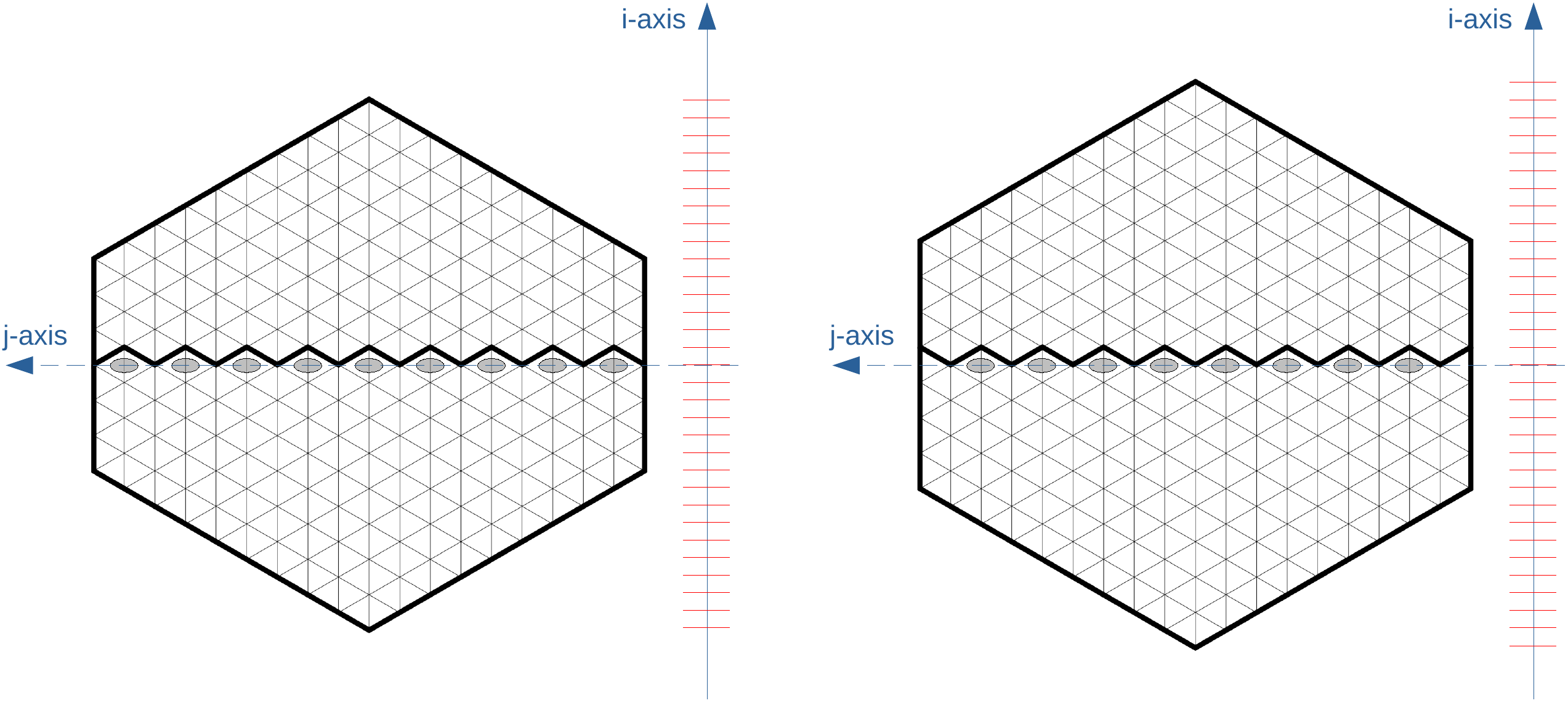}
    \caption{Pictures that show how the Factorization Theorem is applied on the regions $H_{6,9,9}$ and $H_{7,9,9}$. Left picture shows how the Factorization Theorem split the region $H_{6,9,9}$ into two subregions $H_{6,9,9}^{+}$ (bottom) and $H_{6,9,9}^{-}$ (top). The picture on the right present $H_{7,9,9}^{+}$ (bottom) and $H_{7,9,9}^{-}$ (top) that can be obtained from $H_{7,9,9}$ by applying the theorem.}
\end{figure}

We split the region $H_{o}(a,m,m,k:a_1,\ldots,a_k)$ in a similar way (see the right picture in Figure 3.3. Since the midpoint of the vertical side of the deleted triangle of side length $2a_k+1$ is not a lattice point, our zigzag line starts at the lattice point that lies just above the midpoint). Again, $H_{o}^{+}(a,m,m,k:a_1,\ldots,a_k)$ and $H_{o}^{-}(a,m,m,k:a_1,\ldots,a_k)$ are regions on the bottom and the top, and we consider their TGFs $M_{q}(H_{o}^{+}(a,m,m,k:a_1,\ldots,a_k))$ and $M_{q}(H_{o}^{-}(a,m,m,k:a_1,\ldots,a_k))$. As in the previous case, we use the weight of the lozenges inherited from the weight assignment of the lozenges in $H_{o}(a,m,m,k:a_1,\ldots,a_k)$, except the horizontal lozenges in $H_{o}^{+}(a,m,m,k:a_1,\ldots,a_k)$ that are adjacent to the zigzag line (instead of $\frac{q^{0}+q^{0}}{2}=1$, we assign a weight $\frac{1}{2}$ to them). Then, the Matching Factorization Theorem gives the following factorization of $M_{q}(H_{o}(a,m,m,k:a_1,\ldots,a_k))$:

\begin{equation}
\begin{aligned}
    &M_{q}(H_{o}(a,m,m,k:a_1,\ldots,a_k))\\
    &=2^{m}M_{q}(H_{o}^{+}(a,m,m,k:a_1,\ldots,a_k))M_{q}(H_{o}^{-}(a,m,m,k:a_1,\ldots,a_k)).
\end{aligned}
\end{equation}

We again apply the Factorization Theorem to the hexagonal region $H_{x,y,y}$. We divide the region $H_{x,y,y}$ into two subregions using a similar zigzag line and assign weight $\frac{1}{2}$ (instead of $1$) to horizontal lozenges on the symmetry axis (see Figure 3.4). For all other lozenges, we use the weight that were inherited from the weight assignment of the lozenges in $H_{x,y,y}$. The regions at the bottom and the top are indicated by $H_{x,y,y}^{+}$ and $H_{x,y,y}^{-}$, respectively, and their tiling generating functions are denoted by $M_{q}(H_{x,y,y}^{+})$ and $M_{q}(H_{x,y,y}^{-})$. Then $M_{q}(H_{x,y,y})$ has the following factorization:

\begin{equation}
    M_{q}(H_{x,y,y})=2^{y}M_{q}(H_{x,y,y}^{+})M_{q}(H_{x,y,y}^{-}).
\end{equation}

The last theorems we need are two versions of Kuo's graphical condensation method~\cite{Kuo}. This method is an alternative combinatorial interpretation of the Kasteleyn-Percus method~\cite{Ka, Pe}, as Fulmek pointed out in~\cite{Fu1}. To state the theorem, denote a bipartite graph by $G=(V_1, V_2, E)$, where $E$ is the set of edges of the graph and $(V_1, V_2)$ is the partition of the vertex set of the graph $G$ such that every edge in $E$ connects a vertex in $V_1$ and a vertex in $V_2$. Also, for any set of vertices $\{x_1,\ldots,x_n\}$, let $G-\{x_1,\ldots,x_n\}$ be the graph obtained from $G$ by deleting $n$ vertices $x_1,\ldots,x_n$ and all edges adjacent to at least one of $x_1,\ldots,x_n$. \textit{A perfect matching} of a bipartite graph is a subset of edges such that every vertex is incident to precisely one edge. If all edges of a graph are weighted, \textit{the weight of a perfect matching} is defined to be the product of the weights of all edges that constitute the perfect matching. For any weighted graph $G$, \textit{a matching generating function (MGF)} is the sum of weights of all its perfect matchings and is denoted by $M(G)$.

\begin{thm}[\cite{Kuo}, Theorem 5.1]
Let $G=(V_1, V_2, E)$ be a weighted plane bipartite graph in which $|V_1|=|V_2|$. Let vertices $x, y, z,$ and $w$ appear in a cyclic order on a face of $G$. If $x,z \in V_1$ and $y,w \in V_2$, then
\begin{equation}
    M(G)M(G-\{x,y,z,w\})=M(G-\{x,y\})M(G-\{z,w\})+M(G-\{x,w\})M(G-\{y,z\}).
\end{equation}
\end{thm}

\begin{thm}[\cite{Kuo}, Theorem 5.3]
Let $G=(V_1, V_2, E)$ be a weighted plane bipartite graph in which $|V_1|=|V_2|+1$. Let vertices $x, y, z,$ and $w$ appear in a cyclic order on a face of $G$. If $x,y,z \in V_1$ and $w \in V_2$, then
\begin{equation}
    M(G-\{y\})M(G-\{x,z,w\})=M(G-\{x\})M(G-\{y,z,w\})+M(G-\{z\})M(G-\{x,y,w\}).
\end{equation}
\end{thm}

Using a well-known weight-preserving bijection between the set of lozenge tilings of a region and the set of perfect matchings on its dual graph, we apply the two theorems above on the dual graphs of certain regions on a triangular lattice. As a result, we obtain recurrence relations that involve TGFs of lozenges tilings of six related regions on the lattice. Then, we use these recurrence relations to give an inductive proof of our main theorem.

\section{Proof of Theorem 2.2}

The poof is organized as follows: in subsection 4.1, using Ciucu's Matching Factorization Theorem and Lai and Rohatgi's result mentioned in Section 3 (Theorem 3.1), we prove the special case of the theorem when $n=m$. Then, in subsection 4.2, using Kuo's graphical condensation method, we construct three recurrence relations and specialize them to prove another special case when $n=m+1$. We end the proof by proving the general case when $n\geq m+2$ in subsection 4.3 using general forms of the three recurrence relations. When $k=0$, there is nothing to prove because (2.2) and (2.3) become $1=1$. Thus, in the proof, we assume that $k\geq1$.

\subsection{The case $\boldsymbol{n=m}$}
We first prove the case when $n=m$. To prove this case, we need to show that equations (4.1) and (4.2) below hold.

\begin{equation}
\begin{aligned}
    &\frac{M_q(H_{e}(a,m,m,k:a_1,\ldots,a_k))}{M_q(H_{2a+2b_1, m+k, m+k})}\\
    &
\begin{aligned}
    =\prod_{i=0}^{k-1}\Bigg[&\frac{(\langle b_{1}+i+\frac{1}{2}\rangle_{q^2})_{m+k-b_{1}-2i}(\langle a+b_1-b_{i+1}+i+1\langle_{q^2})_{m+k+2b_{i+1}-2i-1}}{(\langle a+b_1+i+\frac{1}{2}\rangle_{q^2})_{m+k-2i}(\langle i+1\rangle_{q^2})_{m+k+b_{i+1}-2i-1}}\\
    &\cdot \frac{(\langle a+m+k+b_1-i\rangle_{q^2})_{b_{i+1}}}{2^{4b_{i+1}}(\langle i+1\rangle_{q^2})_{b_1}(\langle m+k-i-\frac{1}{2}\rangle_{q^2})_{b_{i+1}}(\langle a+b_1+i+1\rangle_{q^2})_{-b_{i+1}}}\\
    &\cdot \prod_{j=1}^{i}\Big[\frac{(\langle b_j-b_{i+1}+i+1-j\rangle_{q^2})^2}{(\langle i+1-j\rangle_{q^2})^2}\cdot(\langle b_{j+1}+i-j+\frac{1}{2}\rangle_{q^2})_{a_j}(\langle b_{j+1}+i-j+1\rangle_{q^2})_{a_j}\Big]\Bigg]
\end{aligned}
\end{aligned}
\end{equation}

and

\begin{equation}
\begin{aligned}
    &\frac{M_q(H_{o}(a,m,m,k:a_1,\ldots,a_k))}{M_q(H_{2a+2b_1+1, m+k, m+k})}\\
    &
\begin{aligned}
    =\prod_{i=0}^{k-1}\Bigg[&\frac{(\langle a+b_1-b_{i+1}+i+1\rangle_{q^2})_{m+k+2b_{i+1}-2i}(\langle b_1+i+\frac{3}{2}\rangle_{q^2})_{m+k-b_1-2i-2}}{(\langle i+1\rangle_{q^2})_{m+k+b_{i+1}-2i-1}(\langle a+b_1+i+\frac{3}{2}\rangle_{q^2})_{m+k-2i-1}}\\
    &\cdot \frac{(\langle b_{i+1}+1\rangle_{q^2})(\langle a+m+k+b_1-i+1\rangle_{q^2})_{b_{i+1}}}{2^{4b_{i+1}+2}(\langle i+2\rangle_{q^2})_{b_1}(\langle m+k-i+\frac{1}{2}\rangle_{q^2})_{b_{i+1}}(\langle a+b_1+i+1\rangle_{q^2})_{-b_{i+1}}}\\
    &
    \begin{aligned}
        \cdot \prod_{j=1}^{i}&\Big[\frac{(\langle b_j-b_{i+1}+i+1-j\rangle_{q^2})^2}{(\langle i+1-j\rangle_{q^2})(\langle b_j+i+1-j\rangle_{q^2})}\frac{(\langle b_j+i+2-j\rangle_{q^2})}{(\langle i+2-j\rangle_{q^2})}\\
        &\cdot (\langle b_{j+1}+i-j+\frac{3}{2}\rangle_{q^2})_{a_j}(\langle b_{j+1}+i-j+2\rangle_{q^2})_{a_j}\Big]\Bigg].
    \end{aligned}
\end{aligned}
\end{aligned}
\end{equation}

First, we show (4.1). The left side of (4.1) can be expressed using telescoping as follows\footnote{When $i=1$, the sequence $(a_1,\ldots,a_{i-1}, b_{i})$ is understood as $(b_{1})$.}\footnote{Here, we need to use the fact that $b_k=\sum_{j=k}^{k}a_j=a_k$.}:

\begin{equation}
\begin{aligned}
    &\frac{M_{q}(H_{e}(a,m,m,k:a_1,\ldots,a_k))}{M_{q}(H_{2a+2b_1, m+k, m+k})}\\
    &
\begin{aligned}
    =&\Bigg[\prod_{i=1}^{k-1}\frac{M_{q}(H_{e}(a,m+k-i-1,m+k-i-1,i+1:a_1,\ldots,a_{i},b_{i+1} ))}{M_{q}(H_{e}(a,m+k-i,m+k-i,i:a_1,\ldots,a_{i-1}, b_{i}))}\Bigg]\\
    &\cdot\frac{M_{q}(H_{e}(a,m+k-1,m+k-1,1:b_1))}{M_{q}(H_{2a+2b_1, m+k, m+k})}.    
\end{aligned}
\end{aligned}
\end{equation}
We investigate each ratio on the right side of (4.3). As mentioned in the previous section, for any $i\in[k]$, by the Matching Factorization Theorem (3.3), $M_{q}(H_{e}(a,m+k-i,m+k-i,i:a_1,\ldots,a_{i-1}, b_{i}))$ can be factorized as follows:

\begin{equation}
\begin{aligned}
    &M_{q}(H_{e}(a,m+k-i,m+k-i,i:a_1,\ldots,a_{i-1}, b_{i}))\\
    &
\begin{aligned}    
    =2^{m+k-i}&M_{q}(H_{e}^{+}(a,m+k-i,m+k-i,i:a_1,\ldots,a_{i-1}, b_{i}))\\
    &\cdot M_{q}(H_{e}^{-}(a,m+k-i,m+k-i,i:a_1,\ldots,a_{i-1}, b_{i})).
\end{aligned}
\end{aligned}
\end{equation}
Similarly, $M_{q}(H_{2a+2b_1, m+k, m+k})$ also has a similar factorization as below:
\begin{equation}
    M_{q}(H_{2a+2b_1, m+k, m+k})=2^{m+k}M_{q}(H_{2a+2b_1, m+k, m+k}^{+})M_{q}(H_{2a+2b_1, m+k, m+k}^{-}).
\end{equation}
Now, we claim that the TGFs of the regions $H_{e}^{+}(a,m+k-i,m+k-i,i:a_1,\ldots,a_{i-1}, b_{i})$, $H_{e}^{-}(a,m+k-i,m+k-i,i:a_1,\ldots,a_{i-1}, b_{i})$, $H_{2a+2b_1, m+k, m+k}^{+}$, and $H_{2a+2b_1, m+k, m+k}^{-}$ are the same as those of certain $R$ regions and $\overline{R}$ regions that we mentioned in the previous section. For any $i\in[k]$, we assert that the following four equations hold.
\begin{equation}
    M_{q}(H_{e}^{+}(a,m+k-i,m+k-i,i:a_1,\ldots,a_{i-1}, b_{i}))=M_{q}(R_{\{b_{i}, b_{i-1}+1, \ldots, b_{1}+i-1\},m+k-i,a}),
\end{equation}
\begin{equation}
\begin{aligned}
    &M_{q}(H_{e}^{-}(a,m+k-i,m+k-i,i:a_1,\ldots,a_{i-1}, b_{i}))\\
    &=M_{q}(\overline{R}_{\{1,2,\ldots,m+k-i-1,m+k-i+b_{i},m+k-i+b_{i-1}+1,\ldots,m+k+b_{1}-1\},0,a}),
\end{aligned}
\end{equation}
\begin{equation}
    M_{q}(H_{2a+2b_1, m+k, m+k}^{+})=M_{q}(R_{\emptyset,m+k,a+b_1-1}),
\end{equation}
and
\begin{equation}
    M_{q}(H_{2a+2b_1, m+k, m+k}^{-})=M_{q}(\overline{R}_{[m+k-1],0,a+b_1}).
\end{equation}
To check the identities above, we first need to observe that one can obtain the identical region (or a mirror image about the $j-$axis) from the corresponding regions after removing some forced lozenges\footnote{A \textit{forced lozenge} on a region $R$ is a lozenge-shaped tile that is contained in all tilings of $R$. In this paper, we indicate forced lozenges by shading them. Suppose that the region $R$ has a forced lozenge, denoted by $L$, which has weight $wt(L)$. Since $L$ is always part of the tilings of $R$, one can check that $M(R-L)=M(R)/wt(L)$ holds, where $M(R)$ and $M(R-L)$ are TGFs of the regions $R$ and $R-L$, respectively. Especially, we have $M(R-L)=M(R)$ when $wt(L)=1$. In this paper, all forced lozenges have a weight $1$, so the TGFs remain unchanged after removing them from regions.}  (see the pictures in the first column of Figure 4.1 that illustrate (4.6) and (4.7)). Furthermore, all removed lozenges are not horizontal; therefore, they have a weight of $1$, so the TGFs are unchanged after removing them. Also, since the weight assignment is symmetric about the $j$-axis, all the regions and their mirror images about $j$-axis have the same TGFs. Therefore, the TGFs of the corresponding regions are the same.

\begin{figure}
    \centering
    \includegraphics[width=16.5cm]{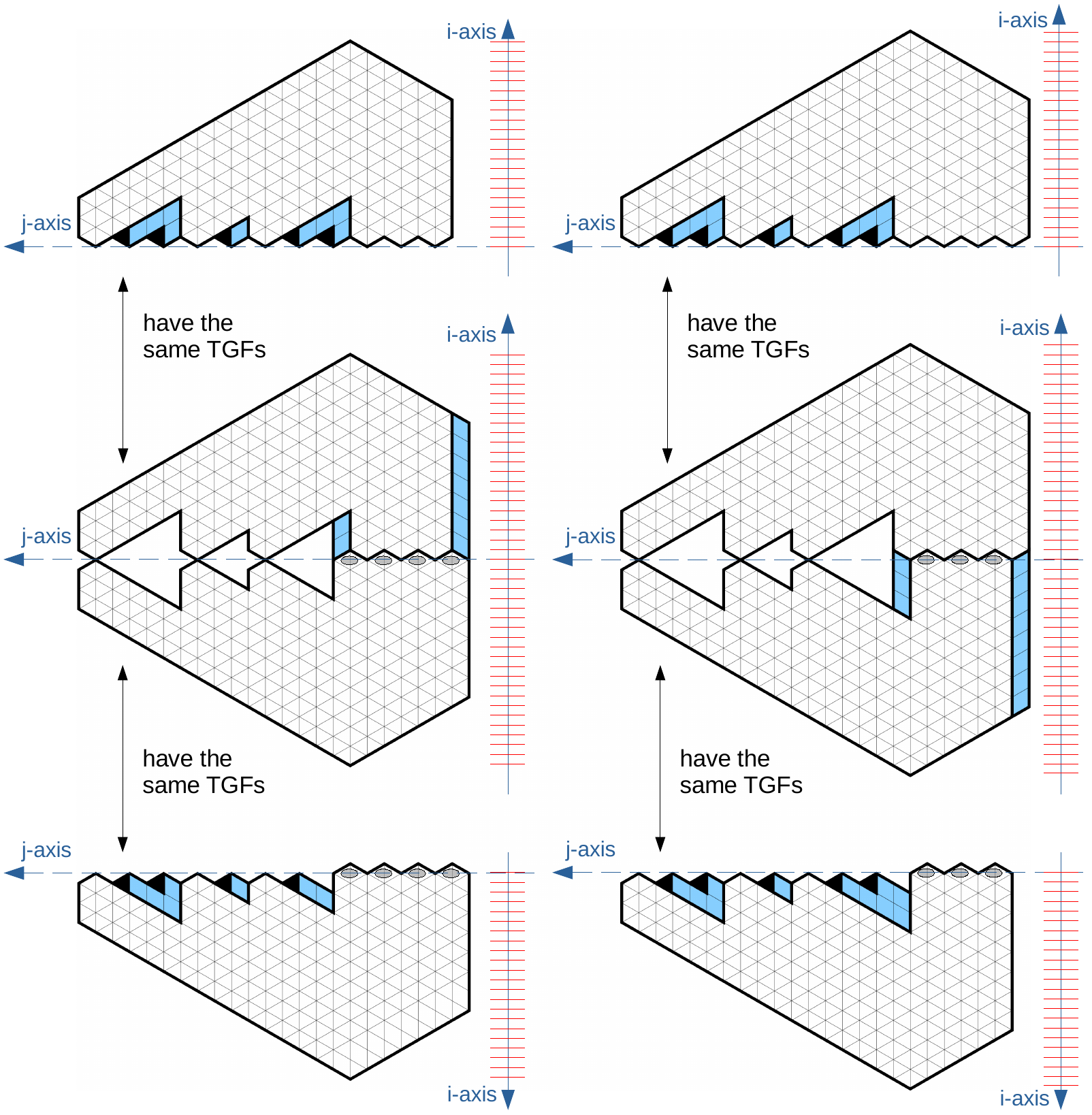}
    \caption{The first column consists of three pictures of $\overline{R}_{\{1,2,3,6,8,11\},0,2}$, $H_{e}(2,4,4,3:2,1,2)$, and a mirror image of $R_{\{2,4,7\},4,2}$ (from top to bottom) and the middle picture is split into $H_{e}^{-}(2,4,4,3:2,1,2)$ and $H_{e}^{+}(2,4,4,3:2,1,2)$ along the zigzag line touching $j$-axis as in Figure 3.3. The second column consists of $\overline{R}_{\{1,2,3,4,7,9,12\},0,2}$, $H_{o}(2,4,4,3:2,1,2)$, and a mirror image of $R_{\{3,5,8\},3,2}$ (from top to bottom) and the middle picture is split into $H_{o}^{-}(2,4,4,3:2,1,2)$ and $H_{o}^{+}(2,4,4,3:2,1,2)$. Forced lozenges are represented by shaded lozenges. One can compare the corresponding regions and check that they have the same TGFs.}

\end{figure}

Using Theorem 3.1, each factor on the right side of (4.3) could be simplified as follows:

\begin{equation}
\begin{aligned}
    &\frac{M_{q}(H_{e}(a,m+k-i-1,m+k-i-1,i+1:a_1,\ldots,a_{i},b_{i+1}))}{M_{q}(H_{e}(a,m+k-i,m+k-i,i:a_1,\ldots,a_{i-1}, b_{i}))}\\
    &
\begin{aligned}
    =&\frac{2^{m+k-i-1}}{2^{m+k-i}}\frac{M_{q}(H_{e}^{+}(a,m+k-i-1,m+k-i-1,i+1:a_1,\ldots,a_{i},b_{i+1}))}{M_{q}(H_{e}^{+}(a,m+k-i,m+k-i,i:a_1,\ldots,a_{i-1}, b_{i}))}\\
    &\cdot\frac{M_{q}(H_{e}^{-}(a,m+k-i-1,m+k-i-1,i+1:a_1,\ldots,a_{i},b_{i+1}))}{M_{q}(H_{e}^{-}(a,m+k-i,m+k-i,i:a_1,\ldots,a_{i-1}, b_{i}))}
\end{aligned}
    \\
    &
\begin{aligned}
    =&\frac{1}{2}\frac{M_{q}(R_{\{b_{i+1}, b_{i}+1, \ldots, b_{1}+i\},m+k-i-1,a})}{M_{q}(R_{\{b_{i}, b_{i-1}+1, \ldots, b_{1}+i-1\},m+k-i,a})}\\
    &\cdot\frac{M_{q}(\overline{R}_{\{1,2,\ldots,m+k-i-2,m+k-i+b_{i+1}-1,m+k-i+b_{i},\ldots,m+k+b_{1}-1\},0,a})}{M_{q}(\overline{R}_{\{1,2,\ldots,m+k-i-1,m+k-i+b_{i},m+k-i+b_{i-1}+1,\ldots,m+k+b_{1}-1\},0,a})}
\end{aligned}    
    \\
    &
\begin{aligned}
    =&\frac{(\langle b_{1}+i+\frac{1}{2}\rangle_{q^2})_{m+k-b_{1}-2i}(\langle a+b_1-b_{i+1}+i+1\rangle_{q^2})_{m+k+2b_{i+1}-2i-1}}{(\langle a+b_1+i+\frac{1}{2}\rangle_{q^2})_{m+k-2i}(\langle i+1\rangle_{q^2})_{m+k+b_{i+1}-2i-1}}\\
    &\cdot\frac{(\langle a+m+k+b_1-i\rangle_{q^2})_{b_{i+1}}}{2^{4b_{i+1}}(\langle i+1\rangle_{q^2})_{b_1}(\langle m+k-i-\frac{1}{2}\rangle_{q^2})_{b_{i+1}}(\langle a+b_1+i+1\rangle_{q^2})_{-b_{i+1}}}\\
    &\cdot\prod_{j=1}^{i}\Big[\frac{(\langle b_j-b_{i+1}+i+1-j\rangle_{q^2})^2}{(\langle i+1-j\rangle_{q^2})^2}\cdot(\langle b_{j+1}+i-j+\frac{1}{2}\rangle_{q^2})_{a_j}(\langle b_{j+1}+i-j+1\rangle_{q^2})_{a_j}\Big]    
\end{aligned}
\end{aligned}
\end{equation}
and

\begin{equation}
\begin{aligned}
    &\frac{M_{q}(H_{e}(a,m+k-1,m+k-1,1:b_1))}{M_{q}(H_{2a+2b_1, m+k, m+k})}\\
    &=\frac{2^{m+k-1}}{2^{m+k}}\frac{M_{q}(H_{e}^{+}(a,m+k-1,m+k-1,1:b_1))}{M_{q}(H_{2a+2b_1, m+k, m+k}^{+})}\frac{M_{q}(H_{e}^{-}(a,m+k-1,m+k-1,1:b_1))}{M_{q}(H_{2a+2b_1, m+k, m+k}^{-})}\\
    &=\frac{1}{2}\frac{M_{q}(R_{\{b_{1}\},m+k-1,a})}{M_{q}(R_{\emptyset,m+k,a+b_1-1})}\frac{M_{q}(\overline{R}_{\{1,2,\ldots,m+k-2,m+k+b_{1}-1\},0,a})}{M_{q}(\overline{R}_{[m+k-1],0,a+b_1})}\\
    &
    \begin{aligned}
        =&\frac{(\langle b_{1}+\frac{1}{2}\rangle_{q^2})_{m+k-b_{1}}(\langle a+1\rangle_{q^2})_{m+k+2b_{1}-1}}{(\langle a+b_1+\frac{1}{2}\rangle_{q^2})_{m+k}(\langle 1\rangle_{q^2})_{m+k+b_{1}-1}}\cdot\frac{(\langle a+m+k+b_1\rangle_{q^2})_{b_{1}}}{2^{4b_{1}}(\langle 1\rangle_{q^2})_{b_1}(\langle m+k-\frac{1}{2}\rangle_{q^2})_{b_{1}}(\langle a+b_1+1\rangle_{q^2})_{-b_{1}}}.
    \end{aligned}
\end{aligned}
\end{equation}
Combining (4.3), (4.10), and (4.11), we get (4.1).

The proof of (4.2) is almost the same as that of (4.1). Instead of (4.6)-(4.9), we use the following identities (see the pictures in the second column of Figure 4.1 that illustrate (4.12) and (4.13)):
\begin{equation}
    M_{q}(H_{o}^{+}(a,m+k-i,m+k-i,i:a_1,\ldots,a_{i-1}, b_{i}))=M_{q}(R_{\{b_{i}+1, b_{i-1}+2, \ldots, b_{1}+i\},m+k-i-1,a}),
\end{equation}
\begin{equation}
\begin{aligned}
    &M_{q}(H_{o}^{-}(a,m+k-i,m+k-i,i:a_1,\ldots,a_{i-1}, b_{i}))\\
    &=M_{q}(\overline{R}_{\{1,2,\ldots,m+k-i,m+k-i+b_{i}+1,m+k-i+b_{i-1}+2,\ldots,m+k+b_{1}\},0,a}),
\end{aligned}
\end{equation}
\begin{equation}
    M_{q}(H_{2a+2b_1+1, m+k, m+k}^{+})=M_{q}(R_{\emptyset,m+k-1,a+b_1}),
\end{equation}
and
\begin{equation}
    M_{q}(H_{2a+2b_1+1, m+k, m+k}^{-})=M_{q}(\overline{R}_{[m+k],0,a+b_1}).
\end{equation}
If one follows the proof of (4.1) (applying the Matching Factorization Theorem (3.4) and using Lai and Rohatgi's result, Theorem 3.1), the proof of (4.2) follows. This completes the proof of (4.1) and (4.2), which are special cases of Theorem 2.2 when $n=m$.

\subsection{The case $\boldsymbol{n=m+1}$}
Now we move on to the proof of the case when $n=m+1$. Our goal is to prove the following two equations:

\begin{equation}
\begin{aligned}
    &\frac{M_{q}(H_{e}(a,m,m+1,k:a_1,\ldots,a_k))}{M_{q}(H_{2a+2b_1, m+k, m+k+1})}\\
    &
\begin{aligned}
    =\prod_{i=0}^{k-1}\Bigg[&\frac{(\langle b_{1}+i+\frac{1}{2}\rangle_{q^2})_{m+k-b_{1}-2i}(\langle a+b_1-b_{i+1}+i+1\rangle_{q^2})_{m+k+2b_{i+1}-2i-1}}{(\langle a+b_1+i+\frac{1}{2}\rangle_{q^2})_{m+k-2i}(\langle i+1\rangle_{q^2})_{m+k+b_{i+1}-2i-1}}\\
    &\cdot \frac{(\langle a+m+k+b_1-i+1\rangle_{q^2})_{b_{i+1}}}{2^{4b_{i+1}}(\langle i+1\rangle_{q^2})_{b_1}(\langle m+k-i+\frac{1}{2}\rangle_{q^2})_{b_{i+1}}(\langle a+b_1+i+1\rangle_{q^2})_{-b_{i+1}}}\\
    &\cdot \prod_{j=1}^{i}\Big[\frac{(\langle b_j-b_{i+1}+i+1-j\rangle_{q^2})^2}{(\langle i+1-j\rangle_{q^2})^2}\cdot(\langle b_{j+1}+i-j+\frac{1}{2}\rangle_{q^2})_{a_j}(\langle b_{j+1}+i-j+1\rangle_{q^2})_{a_j}\Big]\Bigg]
\end{aligned}
\end{aligned}
\end{equation}
and
\begin{equation}
\begin{aligned}
    &\frac{M_{q}(H_{o}(a,m,m+1,k:a_1,\ldots,a_k))}{M_{q}(H_{2a+2b_1+1, m+k, m+k+1})}\\
    &
    \begin{aligned}
    =\prod_{i=0}^{k-1}\Bigg[&\frac{(\langle a+b_1-b_{i+1}+i+1\rangle_{q^2})_{m+k+2b_{i+1}-2i+1}(\langle b_1+i+\frac{3}{2}\rangle_{q^2})_{m+k-b_1-2i-1}}{(\langle i+1\rangle_{q^2})_{m+k+b_{i+1}-2i}(\langle a+b_1+i+\frac{3}{2}\rangle_{q^2})_{m+k-2i}}\\
    &\cdot \frac{(\langle b_{i+1}+1\rangle_{q^2})(\langle a+m+k+b_1-i+1\rangle_{q^2})_{b_{i+1}}}{2^{4b_{i+1}+2}(\langle i+2\rangle_{q^2})_{b_1}(\langle m+k-i+\frac{1}{2}\rangle_{q^2})_{b_{i+1}}(\langle a+b_1+i+1\rangle_{q^2})_{-b_{i+1}}}\\
    &
    \begin{aligned}
        \cdot \prod_{j=1}^{i}&\Big[\frac{(\langle b_j-b_{i+1}+i+1-j\rangle_{q^2})^2}{(\langle i+1-j\rangle_{q^2})(\langle b_j+i+1-j\rangle_{q^2})}\frac{(\langle b_j+i+2-j\rangle_{q^2})}{(\langle i+2-j\rangle_{q^2})}\\
        &\cdot (\langle b_{j+1}+i-j+\frac{3}{2}\rangle_{q^2})_{a_j}(\langle b_{j+1}+i-j+2\rangle_{q^2})_{a_j}\Big]\Bigg].
    \end{aligned}
\end{aligned}
\end{aligned}
\end{equation}

\begin{figure}
    \centering
    \includegraphics[width=15cm]{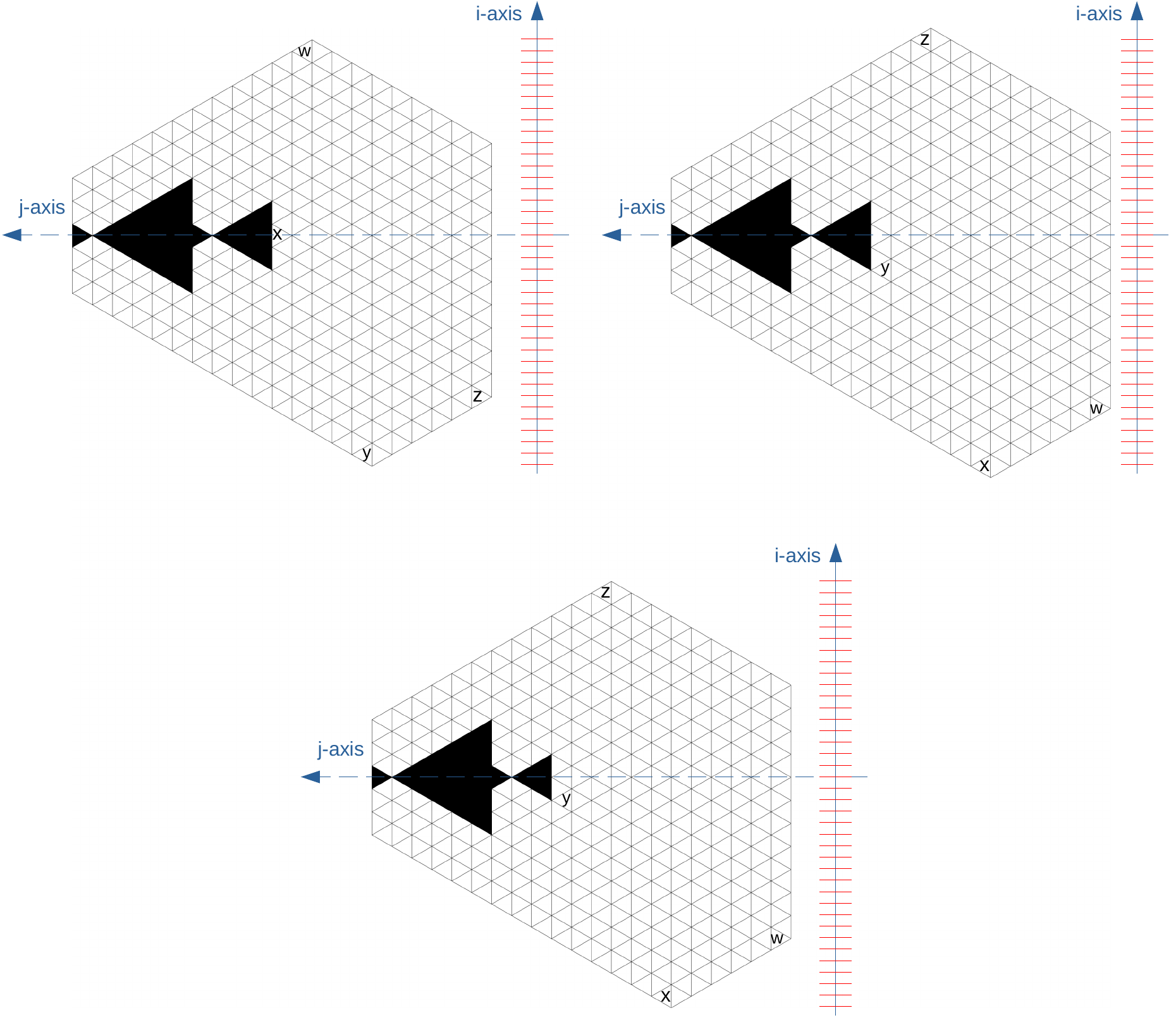}
    \caption{Pictures that show three regions and choices of four unit triangles $x,y,z,$ and $w$ on each region.}
\end{figure}

To achieve the goal, we construct three recurrence relations using Kuo's graphical condensation method to give an inductive proof. First, we consider the region $H_{o}(a,m+1,n,k-1:a_1,\ldots,a_{k-1})$ on the $(i,j)-$plane and choose four unit triangles $x,y,z,$ and $w$ as described in Figure 4.2 (see the upper left picture of the figure). Then we consider the dual graph of the region and four vertices of the dual graph that correspond to four unit triangles $x,y,z,$ and $w$. They satisfy the conditions in Theorem 3.2, so we can apply the theorem on the dual graph. As a result, we obtain an identity that involves MGFs of the dual graph of the region and its subgraphs. However, using a well-known weight-preserving bijection between the set of lozenge tilings of a region and the set of perfect matchings on its dual graph, the identity can be converted into an identity that involves TGFs of $H_{o}(a,m+1,n,k-1:a_1,\ldots,a_{k-1})$ and the same region with some unit triangles removed (see ~\cite{B} for more explanation of this bijection). As can be seen in Figure 4.3, the removal of each unit triangle generates some forced lozenges. However, since all these forced lozenges are not horizontal, they are all weighted by $1$, and removing them from the regions does not change the TGFs\footnote{We use this argument again later when we construct the remaining two recurrence relations.}. Therefore, we obtain the following recurrence relation (see Figure 4.3, which shows how we obtain the six regions):

\begin{figure}
    \centering
    \includegraphics[width=13.5cm]{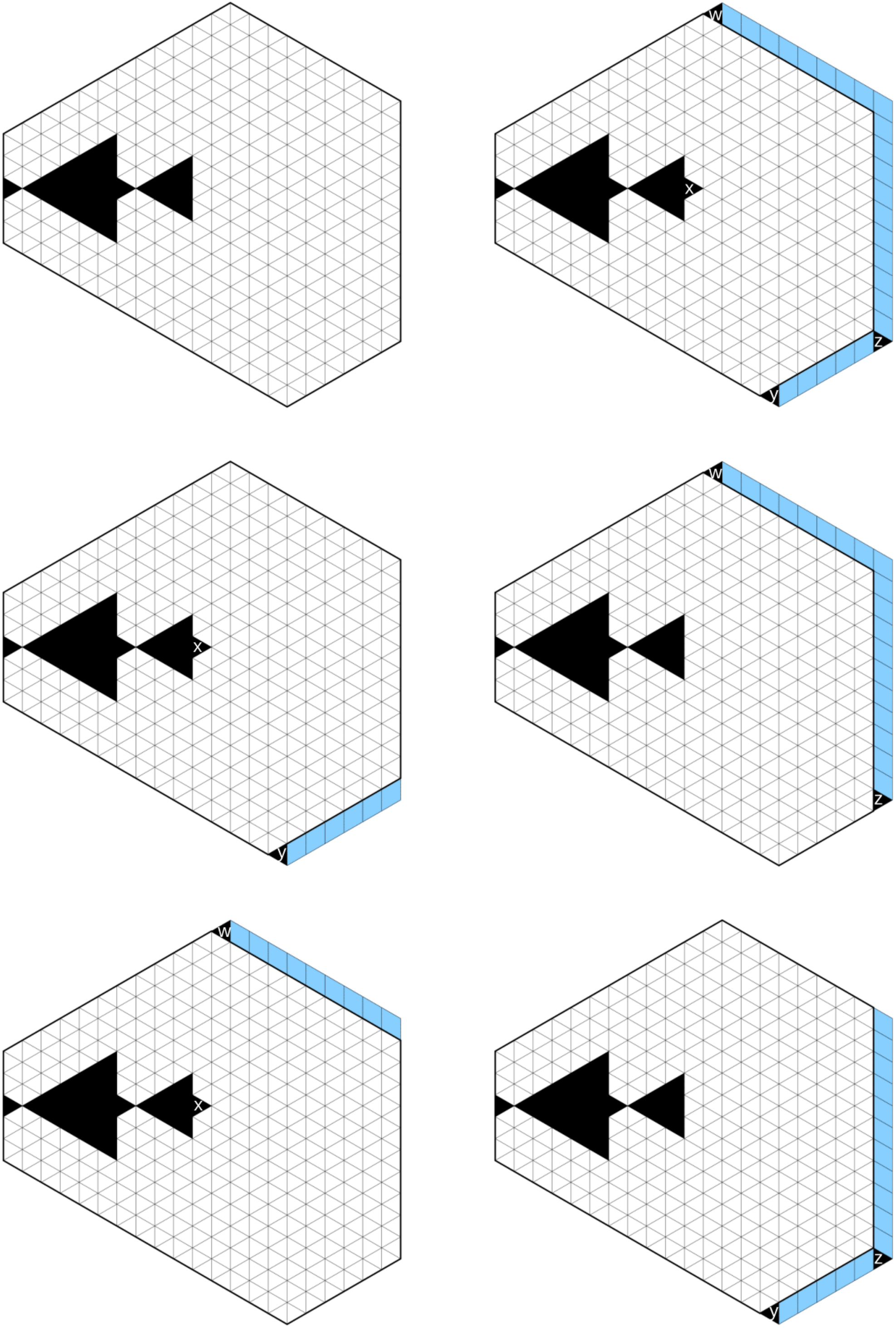}
    \caption{The six regions appearing in the recurrence relation (4.18).}
\end{figure}

\begin{equation}
\begin{aligned}
    &M_{q}(H_{o}(a,m+1,n,k-1:a_1,\ldots,a_{k-1}))M_{q}(H_{e}(a,m,n-1,k:a_1,\ldots,a_{k-1},0))\\
    &
\begin{aligned}
    =&M_{q}(H_{e}(a,m+1,n-1,k:a_1,\ldots,a_{k-1},0))M_{q}(H_{o}(a,m,n,k-1:a_1,\ldots,a_{k-1}))\\
    &+M_{q}(H_{e}(a,m,n,k:a_1,\ldots,a_{k-1},0))M_{q}(H_{o}(a,m+1,n-1,k-1:a_1,\ldots,a_{k-1})).
\end{aligned}
\end{aligned}
\end{equation}

Now we consider the region $H_{e}(a,m+1,n,k:a_1,\ldots,a_k)$ on the $(i,j)-$plane with $a_k\geq 1$. From the region, we replace the removed left-pointing triangle of side length $2a_k$ with that of side length $2a_k-1$. Note that this new region contains one more left-pointing unit triangle than right-pointing unit triangles. From the new region, we choose four unit triangles $x,y,z,$ and $w$ as described in Figure 4.2 (see the picture on the upper right). Note that this time, the dual graph of the region and the four vertices of the dual graph that correspond to the four unit triangles $x,y,z,$ and $w$ satisfy the conditions of Theorem 3.3. If we apply Theorem 3.3 and use the bijection again, after removing forced lozenges of weight 1, we get the following recurrence relation (see Figure 4.4, which shows how we obtain the six regions):

\begin{figure}
    \centering
    \includegraphics[width=13.5cm]{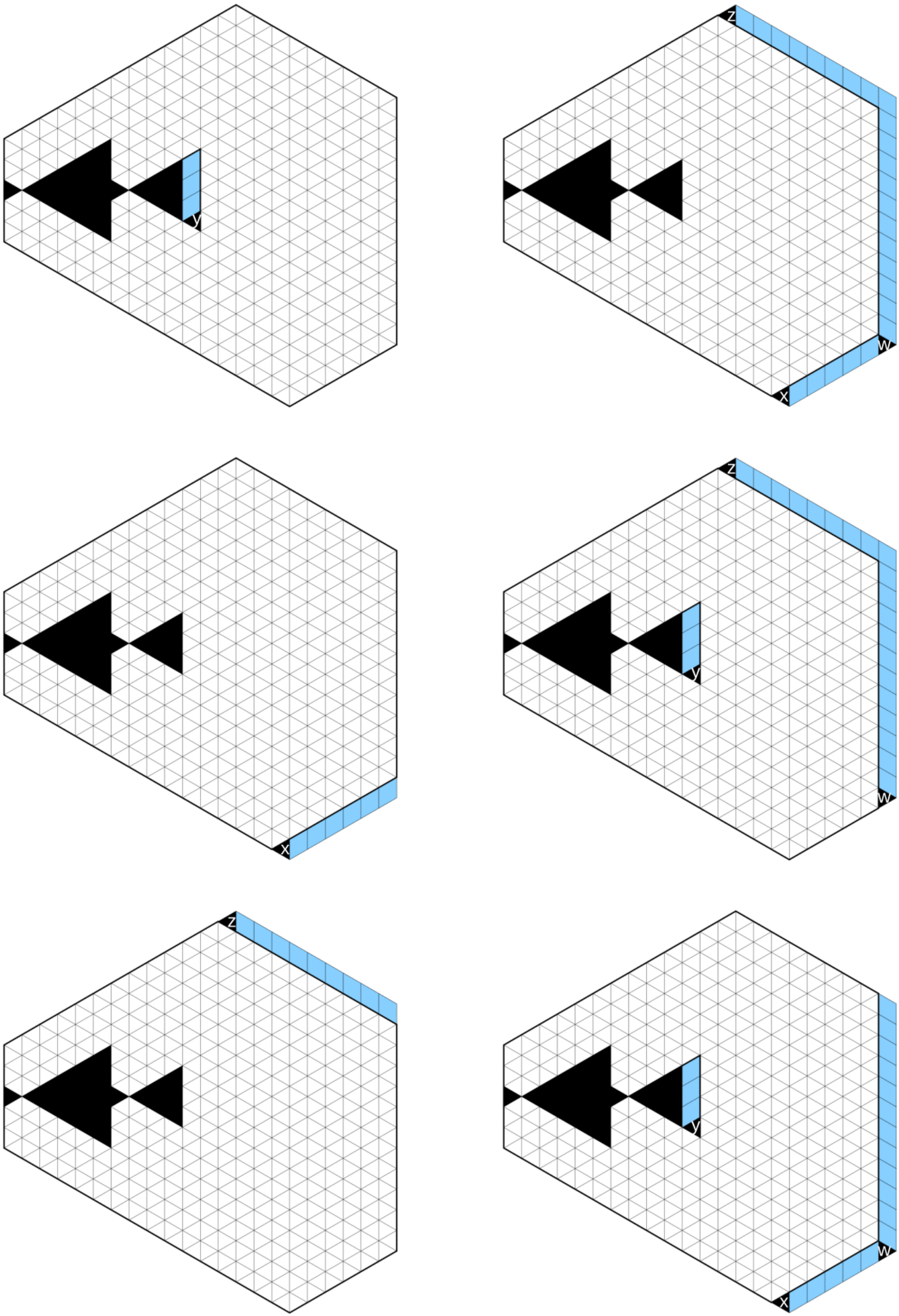}
    \caption{The six regions appearing in the recurrence relation (4.19).}
\end{figure}

\begin{equation}
\begin{aligned}
    &M_{q}(H_{e}(a,m+1,n,k:a_1,\ldots,a_k))M_{q}(H_{o}(a,m+1,n,k:a_1,\ldots,a_k-1))\\
    &
\begin{aligned}
    =&M_{q}(H_{o}(a,m+2,n,k:a_1,\ldots,a_k-1))M_{q}(H_{e}(a,m,n,k:a_1,\ldots,a_k))\\
    &+M_{q}(H_{o}(a,m+1,n+1,k:a_1,\ldots,a_k-1))M_{q}(H_{e}(a,m+1,n-1,k:a_1,\ldots,a_k)).
\end{aligned}
\end{aligned}
\end{equation}

For the last recurrence relation, we consider the region $H_{o}(a,m+1,n,k:a_1,\ldots,a_k)$ on the $(i,j)-$plane (this time, we allow $a_k$ to be $0$). From the region, we replace the removed left-pointing triangle of side length $2a_k+1$ with that of side length $2a_k$. From the new region, we choose four unit triangles $x,y,z,$ and $w$ as described in Figure 4.2 (see the bottom picture in the figure). Using the same bijection used earlier, we can again apply Theorem 3.3 and generate the following recurrence relation (see Figure 4.5, which shows how we obtain the six regions):

\begin{figure}
    \centering
    \includegraphics[width=13.5cm]{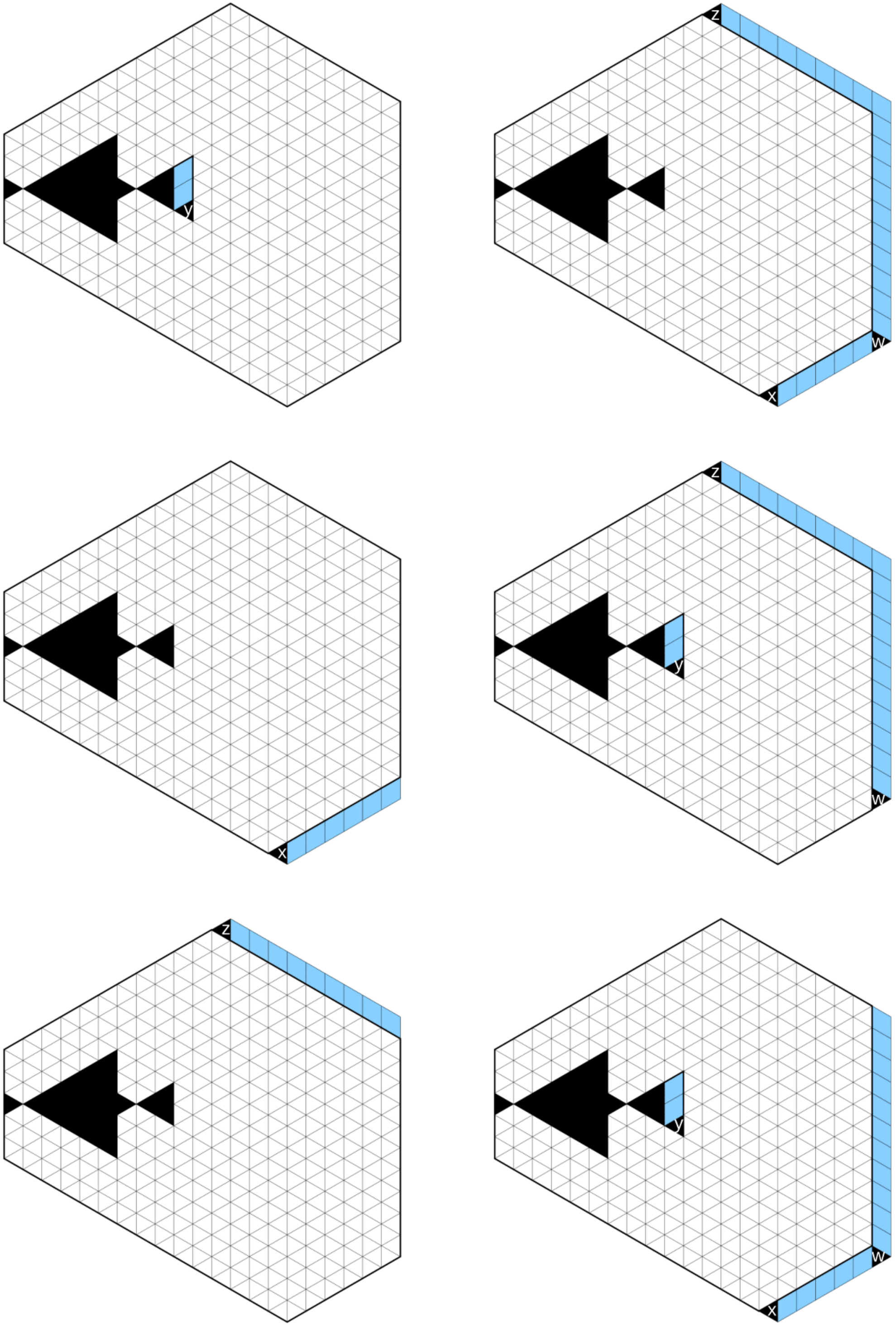}
    \caption{The six regions appearing in the recurrence relation (4.20).}
\end{figure}

\begin{equation}
\begin{aligned}
    &M_{q}(H_{o}(a,m+1,n,k:a_1,\ldots,a_k))M_{q}(H_{e}(a,m+1,n,k:a_1,\ldots,a_k))\\
    &
\begin{aligned}
    =&M_{q}(H_{e}(a,m+2,n,k:a_1,\ldots,a_k))M_{q}(H_{o}(a,m,n,k:a_1,\ldots,a_k))\\
    &+M_{q}(H_{e}(a,m+1,n+1,k:a_1,\ldots,a_k))M_{q}(H_{o}(a,m+1,n-1,k:a_1,\ldots,a_k)).
\end{aligned}
\end{aligned}
\end{equation}

To convert these recurrence relations into suitable forms to prove the case when $n=m+1$, we specialize $n$ in the three recurrence relations above. Note that the original forms of the three recurrence relations (4.18)-(4.20) are used later when we prove the general case for $n\geq m+2$.
We first set $n=m+1$ in (4.18). Then, we get

\begin{equation}
\begin{aligned}
    &M_{q}(H_{o}(a,m+1,m+1,k-1:a_1,\ldots,a_{k-1}))M_{q}(H_{e}(a,m,m,k:a_1,\ldots,a_{k-1},0))\\
    &
\begin{aligned}
    =&M_{q}(H_{e}(a,m+1,m,k:a_1,\ldots,a_{k-1},0))M_{q}(H_{o}(a,m,m+1,k-1:a_1,\ldots,a_{k-1}))\\
    &+M_{q}(H_{e}(a,m,m+1,k:a_1,\ldots,a_{k-1},0))M_{q}(H_{o}(a,m+1,m,k-1:a_1,\ldots,a_{k-1})).
\end{aligned}
\end{aligned}
\end{equation}
Note that the regions $H_{e}(a,m+1,m,k:a_1,\ldots,a_{k-1},0)$ and $H_{e}(a,m,m+1,k:a_1,\ldots,a_{k-1},0)$ are mirror images of each other about the $j-$axis, so they have the same TGFs\footnote{This is because our weight assignment on lozenges is symmetric about the $j$-axis.}. By the same reasoning, the two regions $H_{o}(a,m,m+1,k-1:a_1,\ldots,a_{k-1})$ and $H_{o}(a,m+1,m,k-1:a_1,\ldots,a_{k-1})$ also have the same TGFs. Thus, the identity above can be simplified as follows:

\begin{equation}
\begin{aligned}
    &M_{q}(H_{o}(a,m+1,m+1,k-1:a_1,\ldots,a_{k-1}))M_{q}(H_{e}(a,m,m,k:a_1,\ldots,a_{k-1},0))\\
    &=2M_{q}(H_{e}(a,m,m+1,k:a_1,\ldots,a_{k-1},0))M_{q}(H_{o}(a,m,m+1,k-1:a_1,\ldots,a_{k-1})).
\end{aligned}
\end{equation}

Using the same specialization (setting $n=m+1$), one can obtain the following identities from (4.19) and (4.20):
\begin{equation}
\begin{aligned}
    &M_{q}(H_{e}(a,m+1,m+1,k:a_1,\ldots,a_k))M_{q}(H_{o}(a,m+1,m+1,k:a_1,\ldots,a_k-1))\\
    &=2M_{q}(H_{o}(a,m+1,m+2,k:a_1,\ldots,a_k-1))M_{q}(H_{e}(a,m,m+1,k:a_1,\ldots,a_k))
\end{aligned}
\end{equation}
and
\begin{equation}
\begin{aligned}
    &M_{q}(H_{o}(a,m+1,m+1,k:a_1,\ldots,a_k))M_{q}(H_{e}(a,m+1,m+1,k:a_1,\ldots,a_k))\\
    &=2M_{q}(H_{e}(a,m+1,m+2,k:a_1,\ldots,a_k))M_{q}(H_{o}(a,m,m+1,k:a_1,\ldots,a_k)).
\end{aligned}
\end{equation}

Using these three identities, we now prove the case when $n=m+1$. To give a proof, we first define a concept of \textit{depth} (denoted by $d$) on two families of regions $H_{e}(a,m,n,k:a_1,\ldots,a_k))$ and $H_{o}(a,m,n,k:a_1,\ldots,a_k))$. The depth of the regions $H_{e}(a,m,n,k:a_1,\ldots,a_k))$ and $H_{o}(a,m,n,k:a_1,\ldots,a_k))$ is defined to be $b_1+k=(\sum_{j=1}^{k}a_j)+k$. We prove the case when $n=m+1$ using an induction on the depth $d$ we just defined. When $d=0$, there is nothing to prove because $d=0$ implies $k=0$ and the regions $H_{e}(a,m,m+1,0:\cdot)$ and $H_{o}(a,m,m+1,0:\cdot)$ are the same as $H_{2a,m,m+1}$ and $H_{2a+1,m,m+1}$, respectively. Thus, the identities (4.16) and (4.17) become $1=1$, which verifies the case $d=0$. Now suppose that (4.16) and (4.17) are true for any regions with $d<D$ for some positive integer $D$. Under this assumption, we show that the same holds for the regions with $d=D$. We check this by following the three steps described below:

(Step 1) Show that $M_{q}(H_{e}(a,m,m+1,k:a_1,\ldots,a_{k-1},0))$ is given by (4.16) when the depth of the region is $D$. In the proof of this step, we only use the induction hypothesis.

(Step 2) Show that $M_{q}(H_{e}(a,m,m+1,k:a_1,\ldots,a_{k-1},a_k))$ is given by (4.16) for any $a_k\geq1$ when the depth of the region is $D$. In the proof of this step, we only use the induction hypothesis.

(Step 3) Show that $M_{q}(H_{o}(a,m,m+1,k:a_1,\ldots,a_k))$ is given by (4.17) for any $a_k\geq0$ when the depth of the region is $D$. To prove this step, we use both the induction hypothesis and the results of Steps 1 and Step 2.

First, we check Step 1. We use the recurrence relation (4.22) to verify Step 1. In (4.22), the two terms on the left side are known to be given by (4.1) and (4.2) because $m+1=m+1$ and $m=m$. In addition, the depths of $H_{e}(a,m,m+1,k:a_1,\ldots,a_{k-1},0)$ and $H_{o}(a,m,m+1,k-1:a_1,\ldots,a_{k-1})$ are $D$ and $D-1$, respectively. By the induction hypothesis, we know that $M_{q}(H_{o}(a,m,m+1,k-1:a_1,\ldots,a_{k-1}))$ is given by (4.17). Thus, if we check that (4.1), (4.2), (4.16), (4.17), and (2.1) satisfy the recurrence relation (4.22), we can conclude that $M_{q}(H_{e}(a,m,m+1,k:a_1,\ldots,a_{k-1},0))$ is given by (4.16). Note that (4.22) is equivalent to the following equation:

\begin{equation}
    \frac{1}{2}=\frac{M_{q}(H_{e}(a,m,m+1,k:a_1,\ldots,a_{k-1},0))}{M_{q}(H_{e}(a,m,m,k:a_1,\ldots,a_{k-1},0))}\frac{M_{q}(H_{o}(a,m,m+1,k-1:a_1,\ldots,a_{k-1}))}{M_{q}(H_{o}(a,m+1,m+1,k-1:a_1,\ldots,a_{k-1}))}.
\end{equation}
From (4.1), (4.2), (4.16), (4.17), and (2.1),
\begin{equation}
\begin{aligned}
    &\frac{M_{q}(H_{e}(a,m,m+1,k:a_1,\ldots,a_{k-1},0))}{M_{q}(H_{e}(a,m,m,k:a_1,\ldots,a_{k-1},0))}\\
    &
    \begin{aligned}
        =&\prod_{i=1}^{m+k}\Bigg[\frac{q^{(m+1+k)-i}+q^{i-(m+1+k)}}{2}\Bigg]\cdot\frac{\langle m+k\rangle_{q}!\langle 2a+2b_1+2m+2k\rangle_{q}!}{\langle 2m+2k\rangle_{q}!\langle 2a+2b_1+m+k\rangle_{q}!}\\
    &\cdot\prod_{i=0}^{k-1}\Bigg[\frac{\langle m+k-i-\frac{1}{2}\rangle_{q^2}\langle a+m+k+b_1+b_{i+1}-i\rangle_{q^2}}{\langle b_{i+1}+m+k-i-\frac{1}{2}\rangle_{q^2}\langle a+m+k+b_1-i\rangle_{q^2}}\Bigg]
    \end{aligned}
\end{aligned}
\end{equation}
and
\begin{equation}
\begin{aligned}
    &\frac{M_{q}(H_{o}(a,m,m+1,k-1:a_1,\ldots,a_{k-1}))}{M_{q}(H_{o}(a,m+1,m+1,k-1:a_1,\ldots,a_{k-1}))}\\
    &
    \begin{aligned}
        =&\prod_{i=1}^{m+k}\Bigg[\frac{2}{q^{(m+k)-i}+q^{i-(m+k)}}\Bigg]\cdot\frac{\langle 2m+2k-1\rangle_{q}!\langle 2a+2b_1+m+k\rangle_{q}!}{\langle m+k-1\rangle_{q}!\langle 2a+2b_1+2m+2k\rangle_{q}!}\\
        &\cdot\prod_{i=0}^{k-2}\Bigg[\frac{\langle b_{i+1}+m+k-i-\frac{1}{2}\rangle_{q^2}\langle a+m+k+b_1-i\rangle_{q^2}}{\langle m+k-i-\frac{1}{2}\rangle_{q^2}\langle a+m+k+b_1+b_{i+1}-i\rangle_{q^2}}\Bigg].
    \end{aligned}
\end{aligned}
\end{equation}
Putting (4.26) and (4.27) together and using the fact that $b_k=a_k=0$, we have

\begin{equation}
\begin{aligned}
    &\frac{M_{q}(H_{e}(a,m,m+1,k:a_1,\ldots,a_{k-1},0))}{M_{q}(H_{e}(a,m,m,k:a_1,\ldots,a_{k-1},0))}\frac{M_{q}(H_{o}(a,m,m+1,k-1:a_1,\ldots,a_{k-1}))}{M_{q}(H_{o}(a,m+1,m+1,k-1:a_1,\ldots,a_{k-1}))}\\
    &=\frac{q^{m+k}+q^{-(m+k)}}{2}\cdot\frac{\langle m+k\rangle_{q}}{\langle 2m+2k\rangle_{q}}\\
    &=\frac{1}{2}.
\end{aligned}
\end{equation}
and thus Step 1 is checked.

We move on to Step 2. Here, we use the recurrence relation (4.23). In (4.23), the two terms on the left side are given by (4.1) and (4.2) because $m+1=m+1$. The depths of the two regions on the right, $H_{o}(a,m+1,m+2,k:a_1,\ldots,a_k-1)$ and $H_{e}(a,m,m+1,k:a_1,\ldots,a_k)$, are $D-1$ and $D$, respectively. By the induction hypothesis, we know that $M_{q}(H_{o}(a,m+1,m+2,k:a_1,\ldots,a_k-1))$ is given by (4.17). Therefore, to show that $M_{q}(H_{e}(a,m,m+1,k:a_1,\ldots,a_k))$ is given by (4.16) when $a_k\geq 1$, it suffices to show that (4.1), (4.2), (4.16), (4.17), and (2.1) satisfy the recurrence relations (4.23). Note that (4.23) is equivalent to the following equation:

\begin{equation}
    \frac{1}{2}=\frac{M_{q}(H_{e}(a,m,m+1,k:a_1,\ldots,a_k))}{M_{q}(H_{e}(a,m+1,m+1,k:a_1,\ldots,a_k))}\frac{M_{q}(H_{o}(a,m+1,m+2,k:a_1,\ldots,a_k-1))}{M_{q}(H_{o}(a,m+1,m+1,k:a_1,\ldots,a_k-1))}.
\end{equation}
From (4.1), (4.2), (4.16), (4.17), and (2.1),
\begin{equation}
\begin{aligned}
    &\frac{M_{q}(H_{e}(a,m,m+1,k:a_1,\ldots,a_k))}{M_{q}(H_{e}(a,m+1,m+1,k:a_1,\ldots,a_k))}\\
    &
    \begin{aligned}
        =&\prod_{i=1}^{m+k+1}\Bigg[\frac{2}{q^{(m+k+1)-i}+q^{i-(m+k+1)}}\Bigg]\cdot\frac{\langle 2m+2k+1\rangle_{q}!\langle 2a+2b_1+m+k\rangle_{q}!}{\langle m+k\rangle_{q}!\langle 2a+2b_1+2m+2k+1\rangle_{q}!}\\
        &\cdot\prod_{i=0}^{k-1}\Bigg[\frac{\langle b_{i+1}+m+k-i\rangle_{q^2}\langle a+b_1+m+k-i+\frac{1}{2}\rangle_{q^2}}{\langle m+k-i+\frac{1}{2}\rangle_{q^2}\langle a+b_1+b_{i+1}+m+k-i\rangle_{q^2}}\Bigg]
    \end{aligned}
\end{aligned}
\end{equation}
and
\begin{equation}
\begin{aligned}
    &\frac{M_{q}(H_{o}(a,m+1,m+2,k:a_1,\ldots,a_k-1))}{M_{q}(H_{o}(a,m+1,m+1,k:a_1,\ldots,a_k-1))}\\
    &
    \begin{aligned}
        =&\prod_{i=1}^{m+k+1}\Bigg[\frac{q^{(m+k+2)-i}+q^{i-(m+k+2)}}{2}\Bigg]\cdot\frac{\langle m+k+1\rangle_{q}!\langle 2a+2b_1+2m+2k+1\rangle_{q}!}{\langle 2m+2k+2\rangle_{q}!\langle 2a+2b_1+m+k\rangle_{q}!}\\
        &\cdot\prod_{i=0}^{k-1}\Bigg[\frac{\langle m+k-i+\frac{1}{2}\rangle_{q^2}\langle a+b_1+b_{i+1}+m+k-i\rangle_{q^2}}{\langle b_{i+1}+m+k-i\rangle_{q^2}\langle a+b_1+m+k-i+\frac{1}{2}\rangle_{q^2}}\Bigg].
    \end{aligned}
\end{aligned}
\end{equation}
Putting (4.30) and (4.31) together,
\begin{equation}
\begin{aligned}
    &\frac{M_{q}(H_{e}(a,m,m+1,k:a_1,\ldots,a_k))}{M_{q}(H_{e}(a,m+1,m+1,k:a_1,\ldots,a_k))}\frac{M_{q}(H_{o}(a,m+1,m+2,k:a_1,\ldots,a_k-1))}{M_{q}(H_{o}(a,m+1,m+1,k:a_1,\ldots,a_k-1))}\\
    &=\frac{q^{m+k+1}+q^{-(m+k+1)}}{2}\cdot\frac{\langle m+k+1\rangle_{q}}{\langle 2m+2k+2\rangle_{q}}\\
    &=\frac{1}{2}.
\end{aligned}
\end{equation}
and thus Step 2 is done.

A slightly different argument is needed to check Step 3. As we already mentioned, we use both the induction hypothesis and the results of Step 1 and Step 2 to show that $M_{q}(H_{o}(a,m,m+1,k:a_1,\ldots,a_k))$ is given by (4.17) when its depth is $D$. In this case, we use the recurrence relation (4.24). In (4.24), two terms on the left side is known to be given by (4.1) and (4.2) because $m+1=m+1$. On the right side, the depths of the two regions that appeared there ($H_{e}(a,m+1,m+2,k:a_1,\ldots,a_k)$ and $H_{o}(a,m,m+1,k:a_1,\ldots,a_k)$) are both $D$. However, in Step 1 and Step 2, we showed that $H_{e}(a,m+1,m+2,k:a_1,\ldots,a_k)$ is given by (4.16) under the induction hypothesis (this is because $m+2=(m+1)+1$ and the depth is $D$). Thus, to prove that $M_{q}(H_{o}(a,m,m+1,k:a_1,\ldots,a_k))$ is given by (4.17), we only have to check that (4.1), (4.2), (4.16), (4.17), and (2.1) satisfy (4.24). Note that (4.24) is equivalent to the following equation:
\begin{equation}
    \frac{1}{2}=\frac{M_{q}(H_{e}(a,m+1,m+2,k:a_1,\ldots,a_k))}{M_{q}(H_{e}(a,m+1,m+1,k:a_1,\ldots,a_k))}\frac{M_{q}(H_{o}(a,m,m+1,k:a_1,\ldots,a_k))}{M_{q}(H_{o}(a,m+1,m+1,k:a_1,\ldots,a_k))}.
\end{equation}
From (4.1), (4.2), (4.16), (4.17), and (2.1),
\begin{equation}
\begin{aligned}
    &\frac{M_{q}(H_{e}(a,m+1,m+2,k:a_1,\ldots,a_k))}{M_{q}(H_{e}(a,m+1,m+1,k:a_1,\ldots,a_k))}\\
    &
    \begin{aligned}
        =&\prod_{i=1}^{m+k+1}\Bigg[\frac{q^{(m+k+2)-i}+q^{i-(m+k+2)}}{2}\Bigg]\cdot\frac{\langle m+k+1\rangle_{q}!\langle 2a+2b_1+2m+2k+2\rangle_{q}!}{\langle 2m+2k+2\rangle_{q}!\langle 2a+2b_1+m+k+1\rangle_{q}!}\\
        &\cdot\prod_{i=0}^{k-1}\Bigg[\frac{\langle m+k-i+\frac{1}{2}\rangle_{q^2}\langle a+b_1+b_{i+1}+m+k-i+1\rangle_{q^2}}{\langle b_{i+1}+m+k-i+\frac{1}{2}\rangle_{q^2}\langle a+b_1+m+k-i+1\rangle_{q^2}}\Bigg]
    \end{aligned}
\end{aligned}
\end{equation}
and
\begin{equation}
\begin{aligned}
    &\frac{M_{q}(H_{o}(a,m,m+1,k:a_1,\ldots,a_k))}{M_{q}(H_{o}(a,m+1,m+1,k:a_1,\ldots,a_k))}\\
    &
    \begin{aligned}
        =&\prod_{i=1}^{m+k+1}\Bigg[\frac{2}{q^{(m+k+1)-i}+q^{i-(m+k+1)}}\Bigg]\cdot\frac{\langle 2m+2k+1\rangle_{q}!\langle 2a+2b_1+m+k+1\rangle_{q}!}{\langle m+k\rangle_{q}!\langle 2a+2b_1+2m+2k+2\rangle_{q}!}\\
        &\cdot\prod_{i=0}^{k-1}\Bigg[\frac{\langle b_{i+1}+m+k-i+\frac{1}{2}\rangle_{q^2}\langle a+b_1+m+k-i+1\rangle_{q^2}}{\langle m+k-i+\frac{1}{2}\rangle_{q^2}\langle a+b_1+b_{i+1}+m+k-i+1\rangle_{q^2}}\Bigg].
    \end{aligned}
\end{aligned}
\end{equation}
Putting (4.34) and (4.35) together,
\begin{equation}
\begin{aligned}
    &\frac{M_{q}(H_{e}(a,m+1,m+2,k:a_1,\ldots,a_k))}{M_{q}(H_{e}(a,m+1,m+1,k:a_1,\ldots,a_k))}\frac{M_{q}(H_{o}(a,m,m+1,k:a_1,\ldots,a_k))}{M_{q}(H_{o}(a,m+1,m+1,k:a_1,\ldots,a_k))}\\
    &=\frac{q^{m+k+1}+q^{-(m+k+1)}}{2}\cdot\frac{\langle m+k+1\rangle_{q}}{\langle 2m+2k+2\rangle_{q}}\\
    &=\frac{1}{2}.
\end{aligned}    
\end{equation}
Thus, Step 3 is checked, and this completes the proof of the case when $n=m+1$. The two cases when $n=m$ and $n=m+1$ serve as base cases when we prove the general case using mathematical induction.\\

\subsection{The case $\boldsymbol{n\geq m+2}$}
Finally, we finish the proof by showing that the statement is true for any $m$ and $n$ with $n\geq m$. The proof is based on mathematical induction, and we use induction on the value $n-m$. We call this induction \textit{an outer induction}. The cases when $n-m=0$ and $1$ were checked in subsections 4.1 and 4.2, so it is enough to show (2.2) and (2.3) for $n\geq m+2$. Thus, now we assume that (2.2) and (2.3) are true when $n-m<E$ for some positive integer $E\geq2$. Under this assumption, we must show that (2.2) and (2.3) are still true when $n-m=E$.

Strictly speaking, we have to show that (2.2) and (2.3) hold for the cases when $n-m=E$, while the depth of intrusions can be \textit{arbitrary}. Thus, we show the induction step using another induction on depth $d$. We call this induction \textit{an inner induction}. There is nothing to prove when $n-m=E$ and $d=0$ because (2.2) and (2.3) become $1=1$. Suppose that (2.2) and (2.3) are true when $n-m=E$ and $d<D$ for some positive integer $D$. Under this assumption, we must show that (2.2) and (2.3) hold for the regions with $n-m=E$ and $d=D$. Like in the proof of the case when $n=m+1$, we show this by checking the following three steps:

(Step 1$'$) Show that $M_{q}(H_{e}(a,m,n,k:a_1,\ldots,a_{k-1},0))$ is given by (2.2) when $n-m=E$ and the depth of the region is $D$. In the proof of this step, we use the induction hypothesis of both outer and inner induction.

(Step 2$'$) Show that $M_{q}(H_{e}(a,m,n,k:a_1,\ldots,a_{k-1},a_k))$ is given by (2.2) for any $a_k\geq1$ when $n-m=E$ and the depth of the region is $D$. In the proof of this step, we use the induction hypothesis of both outer and inner induction.

(Step 3$'$) Show that $M_{q}(H_{o}(a,m,n,k:a_1,\ldots,a_k))$ is given by (2.3) for any $a_k\geq0$ when $n-m=E$ and the depth of the region is $D$. In the proof of this step, we use the induction hypothesis of both outer and inner induction and the result of Step 1$'$ and Step 2$'$.

First, we check Step 1$'$. To do this, we analyze the regions that appeared in the recurrence relation (4.18). In (4.18), there are six regions appearing in the recurrence relation, and now we claim that five of them are known to be given by (2.2) and (2.3) by the (either outer or inner) induction hypothesis. $M_{q}(H_{o}(a,m+1,n,k-1:a_1,\ldots,a_{k-1}))$, $M_{q}(H_{e}(a,m,n-1,k:a_1,\ldots,a_{k-1},0))$, $M_{q}(H_{e}(a,m+1,n-1,k:a_1,\ldots,a_{k-1},0))$, and $M_{q}(H_{o}(a,m+1,n-1,k-1:a_1,\ldots,a_{k-1}))$ are given by (2.2) and (2.3). It is due to the induction hypothesis of the outer induction because $n-(m+1)=(n-1)-m=n-m-1=E-1<E$ and $(n-1)-(m+1)=n-m-2=E-2<E$. Furthermore, $M_{q}(H_{o}(a,m,n,k-1:a_1,\ldots,a_{k-1}))$ is given by (2.3). This time, it is due to the induction hypothesis of the inner induction because $n-m=E$ and the depth of the region is $D-1$, which is less than $D$. Thus, to show that $M_{q}(H_{e}(a,m,n,k:a_1,\ldots,a_{k-1},0))$ is given by (2.2) when $n-m=E$ and the depth of the region is $D$, it is enough to check that (2.2), (2.3), and (2.1) satisfy the recurrence relation (4.18). Note that (4.18) is equivalent to the following identity:
\begin{equation}
\begin{aligned}
    &\frac{M_{q}(H_{e}(a,m+1,n-1,k:a_1,\ldots,a_{k-1},0))}{M_{q}(H_{e}(a,m,n-1,k:a_1,\ldots,a_{k-1},0))}\cdot\frac{M_{q}(H_{o}(a,m,n,k-1:a_1,\ldots,a_{k-1}))}{M_{q}(H_{o}(a,m+1,n,k-1:a_1,\ldots,a_{k-1}))}\\
    &+\frac{M_{q}(H_{e}(a,m,n,k:a_1,\ldots,a_{k-1},0))}{M_{q}(H_{e}(a,m,n-1,k:a_1,\ldots,a_{k-1},0))}\cdot\frac{M_{q}(H_{o}(a,m+1,n-1,k-1:a_1,\ldots,a_{k-1}))}{M_{q}(H_{o}(a,m+1,n,k-1:a_1,\ldots,a_{k-1}))}\\
    &=1.
\end{aligned}    
\end{equation}
These four ratios that appeared in Equation (4.37) above can be simplified using (2.2), (2.3), and (2.1). After the simplification, one gets
\begin{equation}
\begin{aligned}
    &\frac{M_{q}(H_{e}(a,m+1,n-1,k:a_1,\ldots,a_{k-1},0))}{M_{q}(H_{e}(a,m,n-1,k:a_1,\ldots,a_{k-1},0))}\\
    &
    \begin{aligned}
        =&\prod_{i=1}^{n+k-1}\Bigg[\frac{q^{(m+k+1)-i}+q^{i-(m+k+1)}}{2}\Bigg]\cdot\frac{\langle m+k\rangle_{q}!\langle 2a+2b_1+m+n+2k-1\rangle_{q}!}{\langle m+n+2k-1\rangle_{q}!\langle 2a+2b_1+m+k\rangle_{q}!}\\
        &\cdot\prod_{i=0}^{k-1}\Bigg[\frac{\langle a+b_1+b_{i+1}+m+k-i\rangle_{q^2}\langle \frac{1}{2} m+\frac{1}{2}n+k-i-\frac{1}{2}\rangle_{q^2}\langle \frac{1}{2}m+\frac{1}{2}n+k-i-1\rangle_{q^2}}{\langle m+k-i\rangle_{q^2}\langle a+b_1+\frac{1}{2}m+\frac{1}{2}n+k-i-\frac{1}{2}\rangle_{q^2}\langle b_{i+1}+\frac{1}{2}m+\frac{1}{2}n+k-i-1\rangle_{q^2}}\Bigg],
    \end{aligned}
\end{aligned}
\end{equation}
\begin{equation}
\begin{aligned}
    &\frac{M_{q}(H_{o}(a,m,n,k-1:a_1,\ldots,a_{k-1}))}{M_{q}(H_{o}(a,m+1,n,k-1:a_1,\ldots,a_{k-1}))}\\
    &
    \begin{aligned}
        =&\prod_{i=1}^{n+k-1}\Bigg[\frac{2}{q^{(m+k)-i}+q^{i-(m+k)}}\Bigg]\cdot\frac{\langle m+n+2k-2\rangle_{q}!\langle 2a+2b_1+m+k\rangle_{q}!}{\langle m+k-1\rangle_{q}!\langle 2a+2b_1+m+n+2k-1\rangle_{q}!}\\
        &\cdot\prod_{i=0}^{k-2}\Bigg[\frac{\langle m+k-i-1\rangle_{q^2}\langle a+b_1+\frac{1}{2}m+\frac{1}{2}n+k-i-\frac{1}{2}\rangle_{q^2}\langle b_{i+1}+\frac{1}{2}m+\frac{1}{2}n+k-i-1\rangle_{q^2}}{\langle a+b_1+b_{i+1}+m+k-i\rangle_{q^2}\langle \frac{1}{2}m+\frac{1}{2}n+k-i-\frac{3}{2}\rangle_{q^2}\langle \frac{1}{2}m+\frac{1}{2}n+k-i-1\rangle_{q^2}}\Bigg],
    \end{aligned}
\end{aligned}
\end{equation}
\begin{equation}
\begin{aligned}
    &\frac{M_{q}(H_{e}(a,m,n,k:a_1,\ldots,a_{k-1},0))}{M_{q}(H_{e}(a,m,n-1,k:a_1,\ldots,a_{k-1},0))}\\
    &
    \begin{aligned}
        =&\prod_{i=1}^{m+k}\Bigg[\frac{q^{(n+k)-i}+q^{i-(n+k)}}{2}\Bigg]\cdot\frac{\langle n+k-1\rangle_{q}!\langle 2a+2b_1+m+n+2k-1\rangle_{q}!}{\langle m+n+2k-1\rangle_{q}!\langle 2a+2b_1+n+k-1\rangle_{q}!}\\
        &\cdot\prod_{i=0}^{k-1}\Bigg[\frac{\langle a+b_1+b_{i+1}+n+k-i-1\rangle_{q^2}\langle \frac{1}{2}m+\frac{1}{2}n+k-i-\frac{1}{2}\rangle_{q^2}\langle \frac{1}{2}m+\frac{1}{2}n+k-i-1\rangle_{q^2}}{\langle n+k-i-1\rangle_{q^2}\langle a+b_1+\frac{1}{2}m+\frac{1}{2}n+k-i-\frac{1}{2}\rangle_{q^2}\langle b_{i+1}+\frac{1}{2}m+\frac{1}{2}n+k-i-1\rangle_{q^2}}\Bigg],
    \end{aligned}
\end{aligned}
\end{equation}
and
\begin{equation}
\begin{aligned}
    &\frac{M_{q}(H_{o}(a,m+1,n-1,k-1:a_1,\ldots,a_{k-1}))}{M_{q}(H_{o}(a,m+1,n,k-1:a_1,\ldots,a_{k-1}))}\\
    &
    \begin{aligned}
        =&\prod_{i=1}^{m+k}\Bigg[\frac{2}{q^{(n+k-1)-i}+q^{i-(n+k-1)}}\Bigg]\cdot\frac{\langle m+n+2k-2\rangle_{q}!\langle 2a+2b_1+n+k-1\rangle_{q}!}{\langle n+k-2\rangle_{q}!\langle 2a+2b_1+m+n+2k-1\rangle_{q}!}\\
        &\cdot\prod_{i=0}^{k-2}\Bigg[\frac{\langle n+k-i-2\rangle_{q^2}\langle a+b_1+\frac{1}{2}m+\frac{1}{2}n+k-i-\frac{1}{2}\rangle_{q^2}\langle b_{i+1}+\frac{1}{2}m+\frac{1}{2}n+k-i-1\rangle_{q^2}}{\langle a+b_1+b_{i+1}+n+k-i-1\rangle_{q^2}\langle \frac{1}{2}m+\frac{1}{2}n+k-i-\frac{3}{2}\rangle_{q^2}\langle \frac{1}{2}m+\frac{1}{2}n+k-i-1\rangle_{q^2}}\Bigg].
    \end{aligned}
\end{aligned}
\end{equation}
Putting (4.38)-(4.41) together in (4.37) and using $b_k=a_k=0$, one can check that 
\begin{equation}
\begin{aligned}
    &\frac{M_{q}(H_{e}(a,m+1,n-1,k:a_1,\ldots,a_{k-1},0))}{M_{q}(H_{e}(a,m,n-1,k:a_1,\ldots,a_{k-1},0))}\cdot\frac{M_{q}(H_{o}(a,m,n,k-1:a_1,\ldots,a_{k-1}))}{M_{q}(H_{o}(a,m+1,n,k-1:a_1,\ldots,a_{k-1}))}\\
    &+\frac{M_{q}(H_{e}(a,m,n,k:a_1,\ldots,a_{k-1},0))}{M_{q}(H_{e}(a,m,n-1,k:a_1,\ldots,a_{k-1},0))}\cdot\frac{M_{q}(H_{o}(a,m+1,n-1,k-1:a_1,\ldots,a_{k-1}))}{M_{q}(H_{o}(a,m+1,n,k-1:a_1,\ldots,a_{k-1}))}\\
    &
    \begin{aligned}
        =&\frac{q^{m+k}+q^{-(m+k)}}{q^{n-m-1}+q^{-(n-m-1)}}\cdot\frac{\langle m+k\rangle_{q}}{\langle m+n+2k-1\rangle_{q}}\cdot\frac{\langle a+b_1+m+1\rangle_{q^2}\langle \frac{1}{2}m+\frac{1}{2}n+k-\frac{1}{2}\rangle_{q^2}}{\langle m+k\rangle_{q^2}\langle a+b_1+\frac{1}{2}m+\frac{1}{2}n+\frac{1}{2}\rangle_{q^2}}\\
        &+\frac{q^{n+k-1}+q^{-(n+k-1)}}{q^{n-m-1}+q^{-(n-m-1)}}\cdot\frac{\langle n+k-1\rangle_{q}}{\langle m+n+2k-1\rangle_{q}}\cdot\frac{\langle a+b_1+n\rangle_{q^2}\langle \frac{1}{2}m+\frac{1}{2}n+k-\frac{1}{2}\rangle_{q^2}}{\langle n+k-1\rangle_{q^2}\langle a+b_1+\frac{1}{2}m+\frac{1}{2}n+\frac{1}{2}\rangle_{q^2}}
    \end{aligned}
    \\
    &=\frac{\langle a+b_1+m+1\rangle_{q^2}+\langle a+b_1+n\rangle_{q^2}}{(q^{n-m-1}+q^{-(n-m-1)})\langle a+b_1+\frac{1}{2}m+\frac{1}{2}n+\frac{1}{2}\rangle_{q^2}}\\
    &=1.
\end{aligned}    
\end{equation}
Thus, we just check that $M_{q}(H_{e}(a,m,n,k:a_1,\ldots,a_{k-1},0))$ is given by (2.2) when $n-m=E$ and the depth of the region is $D$. It verifies Step 1$'$.

We move on to the verification of Step 2$'$. In this step, we use the recurrence relation (4.19). Among the six terms in (4.19), we claim that five of them are known to be given by (2.2) and (2.3) by the induction hypotheses. Indeed, $M_{q}(H_{e}(a,m+1,n,k:a_1,\ldots,a_k))$, $M_{q}(H_{o}(a,m+1,n,k:a_1,\ldots,a_k-1))$, $M_{q}(H_{o}(a,m+2,n,k:a_1,\ldots,a_k-1))$, and $M_{q}(H_{e}(a,m+1,n-1,k:a_1,\ldots,a_k))$ are given by (2.2) and (2.3). It is due to the induction hypothesis of the outer induction because $n-(m+1)=n-m-1=E-1<E$ and $n-(m+2)=(n-1)-(m+1)=n-m-2=E-2<E$. In addition, $M_{q}(H_{o}(a,m+1,n+1,k:a_1,\ldots,a_k-1))$ is given by (2.3). It is due to the induction hypothesis of the inner induction because $(n+1)-(m+1)=n-m=E$ and the depth of the region is $D-1<D$. Therefore, to show that $M_{q}(H_{e}(a,m,n,k:a_1,\ldots,a_{k-1},a_k))$ with $a_k\geq1$ is given by (2.2), it is enough to check that (2.2), (2.3), and (2.1) satisfy (4.19). Note that (4.19) is equivalent to the following identity:
\begin{equation}
\begin{aligned}
    &\frac{M_{q}(H_{o}(a,m+2,n,k:a_1,\ldots,a_k-1))}{M_{q}(H_{o}(a,m+1,n,k:a_1,\ldots,a_k-1))}\cdot\frac{M_{q}(H_{e}(a,m,n,k:a_1,\ldots,a_k))}{M_{q}(H_{e}(a,m+1,n,k:a_1,\ldots,a_k))}\\
    &+\frac{M_{q}(H_{o}(a,m+1,n+1,k:a_1,\ldots,a_k-1))}{M_{q}(H_{o}(a,m+1,n,k:a_1,\ldots,a_k-1))}\frac{M_{q}(H_{e}(a,m+1,n-1,k:a_1,\ldots,a_k))}{M_{q}(H_{e}(a,m+1,n,k:a_1,\ldots,a_k))}\\
    &=1.
\end{aligned}
\end{equation}
As before, we can use (2.2), (2.3), and (2.1) to obtain the following identities:
\begin{equation}
\begin{aligned}
    &\frac{M_{q}(H_{o}(a,m+2,n,k:a_1,\ldots,a_k-1))}{M_{q}(H_{o}(a,m+1,n,k:a_1,\ldots,a_k-1))}\\
    &
    \begin{aligned}
        =&\prod_{i=1}^{n+k}\Bigg[\frac{q^{(m+k+2)-i}+q^{i-(m+k+2)}}{2}\Bigg]\cdot\frac{\langle m+k+1\rangle_{q}!\langle 2a+2b_1+m+n+2k\rangle_{q}!}{\langle m+n+2k+1\rangle_{q}!\langle 2a+2b_1+m+k\rangle_{q}!}\\
        &\cdot\prod_{i=0}^{k-1}\Bigg[\frac{\langle a+b_1+b_{i+1}+m+k-i\rangle_{q^2}\langle \frac{1}{2}m+\frac{1}{2}n+k-i\rangle_{q^2}\langle \frac{1}{2}m+\frac{1}{2}n+k-i+\frac{1}{2}\rangle_{q^2}}{\langle m+k-i+1\rangle_{q^2}\langle a+b_1+\frac{1}{2}m+\frac{1}{2}n+k-i\rangle_{q^2}\langle b_{i+1}+\frac{1}{2}m+\frac{1}{2}n+k-i-\frac{1}{2}\rangle_{q^2}}\Bigg],        
    \end{aligned}
\end{aligned}
\end{equation}
\begin{equation}
\begin{aligned}
    &\frac{M_{q}(H_{e}(a,m,n,k:a_1,\ldots,a_k))}{M_{q}(H_{e}(a,m+1,n,k:a_1,\ldots,a_k))}\\
    &
    \begin{aligned}
        =&\prod_{i=1}^{n+k}\Bigg[\frac{2}{q^{(m+k+1)-i}+q^{i-(m+k+1)}}\Bigg]\cdot\frac{\langle m+n+2k\rangle_{q}!\langle 2a+2b_1+m+k\rangle_{q}!}{\langle m+k\rangle_{q}!\langle 2a+2b_1+m+n+2k\rangle_{q}!}\\
        &\cdot\prod_{i=0}^{k-1}\Bigg[\frac{\langle m+k-i\rangle_{q^2}\langle a+b_1+\frac{1}{2}m+\frac{1}{2}n+k-i\rangle_{q^2}\langle b_{i+1}+\frac{1}{2}m+\frac{1}{2}n+k-i-\frac{1}{2}\rangle_{q^2}}{\langle a+b_1+b_{i+1}+m+k-i\rangle_{q^2}\langle \frac{1}{2}m+\frac{1}{2}n+k-i\rangle_{q^2}\langle \frac{1}{2}m+\frac{1}{2}n+k-i-\frac{1}{2}\rangle_{q^2}}\Bigg],    
    \end{aligned}
\end{aligned}
\end{equation}
\begin{equation}
\begin{aligned}
    &\frac{M_{q}(H_{o}(a,m+1,n+1,k:a_1,\ldots,a_k-1))}{M_{q}(H_{o}(a,m+1,n,k:a_1,\ldots,a_k-1))}\\
    &
    \begin{aligned}
        =&\prod_{i=1}^{m+k+1}\Bigg[\frac{q^{(n+k+1)-i}+q^{i-(n+k+1)}}{2}\Bigg]\cdot\frac{\langle n+k\rangle_{q}!\langle 2a+2b_1+m+n+2k\rangle_{q}!}{\langle m+n+2k+1\rangle_{q}!\langle 2a+2b_1+n+k-1\rangle_{q}!}\\
        &\cdot\prod_{i=0}^{k-1}\Bigg[\frac{\langle a+b_1+b_{i+1}+n+k-i-1\rangle_{q^2}\langle \frac{1}{2}m+\frac{1}{2}n+k-i\rangle_{q^2}\langle \frac{1}{2}m+\frac{1}{2}n+k-i+\frac{1}{2}\rangle_{q^2}}{\langle n+k-i\rangle_{q^2}\langle a+b_1+\frac{1}{2}m+\frac{1}{2}n+k-i\rangle_{q^2}\langle b_{i+1}+\frac{1}{2}m+\frac{1}{2}n+k-i-\frac{1}{2}\rangle_{q^2}}\Bigg],    
    \end{aligned}
\end{aligned}
\end{equation}
and
\begin{equation}
\begin{aligned}
    &\frac{M_{q}(H_{e}(a,m+1,n-1,k:a_1,\ldots,a_k))}{M_{q}(H_{e}(a,m+1,n,k:a_1,\ldots,a_k))}\\
    &
    \begin{aligned}
        =&\prod_{i=1}^{m+k+1}\Bigg[\frac{2}{q^{(n+k)-i}+q^{i-(n+k)}}\Bigg]\cdot\frac{\langle m+n+2k\rangle_{q}!\langle 2a+2b_1+n+k-1\rangle_{q}!}{\langle n+k-1\rangle_{q}!\langle 2a+2b_1+m+n+2k\rangle_{q}!}\\
        &\cdot\prod_{i=0}^{k-1}\Bigg[\frac{\langle n+k-i-1\rangle_{q^2}\langle a+b_1+\frac{1}{2}m+\frac{1}{2}n+k-i\rangle_{q^2}\langle b_{i+1}+\frac{1}{2}m+\frac{1}{2}n+k-i-\frac{1}{2}\rangle_{q^2}}{\langle a+b_1+b_{i+1}+n+k-i-1\rangle_{q^2}\langle \frac{1}{2}m+\frac{1}{2}n+k-i\rangle_{q^2}\langle \frac{1}{2}m+\frac{1}{2}n+k-i-\frac{1}{2}\rangle_{q^2}}\Bigg].    
    \end{aligned}
\end{aligned}
\end{equation}
Putting (4.44)-(4.47) together in (4.43),
\begin{equation}
\begin{aligned}
    &\frac{M_{q}(H_{o}(a,m+2,n,k:a_1,\ldots,a_k-1))}{M_{q}(H_{o}(a,m+1,n,k:a_1,\ldots,a_k-1))}\cdot\frac{M_{q}(H_{e}(a,m,n,k:a_1,\ldots,a_k))}{M_{q}(H_{e}(a,m+1,n,k:a_1,\ldots,a_k))}\\
    &+\frac{M_{q}(H_{o}(a,m+1,n+1,k:a_1,\ldots,a_k-1))}{M_{q}(H_{o}(a,m+1,n,k:a_1,\ldots,a_k-1))}\cdot\frac{M_{q}(H_{e}(a,m+1,n-1,k:a_1,\ldots,a_k))}{M_{q}(H_{e}(a,m+1,n,k:a_1,\ldots,a_k))}\\
    &
    \begin{aligned}
        =&\frac{q^{m+k+1}+q^{-(m+k+1)}}{q^{n-m-1}+q^{-(n-m-1)}}\cdot\frac{\langle m+k+1\rangle_{q}}{\langle m+n+2k+1\rangle_{q}}\cdot\frac{\langle m+1\rangle_{q^2}\langle \frac{1}{2}m+\frac{1}{2}n+k+\frac{1}{2}\rangle_{q^2}}{\langle m+k+1\rangle_{q^2}\langle \frac{1}{2}m+\frac{1}{2}n+\frac{1}{2}\rangle_{q^2}}\\
        &+\frac{q^{n+k}+q^{-(n+k)}}{q^{n-m-1}+q^{-(n-m-1)}}\cdot\frac{\langle n+k\rangle_{q}}{\langle m+n+2k+1\rangle_{q}}\cdot\frac{\langle n\rangle_{q^2}\langle \frac{1}{2}m+\frac{1}{2}n+k+\frac{1}{2}\rangle_{q^2}}{\langle n+k\rangle_{q^2}\langle \frac{1}{2}m+\frac{1}{2}n+\frac{1}{2}\rangle_{q^2}}\\
    \end{aligned}
    \\
    &=\frac{\langle m+1\rangle_{q^2}+\langle n\rangle_{q^2}}{(q^{n-m-1}+q^{-(n-m-1)})\langle \frac{1}{2}m+\frac{1}{2}n+\frac{1}{2}\rangle_{q^2}}\\
    &=1.
\end{aligned}
\end{equation}
and this verifies Step 2$'$.

Lastly, we use (4.20) to check Step 3$'$. Among the six terms in (4.20), we claim that four of them are known to be given by (2.2) and (2.3) by the induction hypotheses, and one of them is known to be given by (2.2) due to the results of Step 1$'$ and Step 2$'$. In fact, $M_{q}(H_{o}(a,m+1,n,k:a_1,\ldots,a_k))$, $M_{q}(H_{e}(a,m+1,n,k:a_1,\ldots,a_k))$, $M_{q}(H_{e}(a,m+2,n,k:a_1,\ldots,a_k))$, and $M_{q}(H_{o}(a,m+1,n-1,k:a_1,\ldots,a_k))$ are given by (2.2) and (2.3). This is due to the induction hypothesis of the outer induction because $n-(m+1)=n-m-1=E-1<E$ and $n-(m+2)=(n-1)-(m+1)=n-m-2=E-2<E$. Also, $M_{q}(H_{e}(a,m+1,n+1,k:a_1,\ldots,a_k))$ is given by (2.2). It is due to the results of Step 1$'$ and Step 2$'$ because $(n+1)-(m+1)=n-m=E$ and the depth of the region is $D$. Thus, to show that $M_{q}(H_{o}(a,m,n,k:a_1,\ldots,a_k))$ is given by (2.3), we have to check that (2.2), (2.3), and (2.1) satisfy the recurrence relation (4.20). Note that the recurrence relation (4.20) is equivalent to the following identity:
\begin{equation}
\begin{aligned}
    &\frac{M_{q}(H_{e}(a,m+2,n,k:a_1,\ldots,a_k))}{M_{q}(H_{e}(a,m+1,n,k:a_1,\ldots,a_k))}\cdot\frac{M_{q}(H_{o}(a,m,n,k:a_1,\ldots,a_k))}{M_{q}(H_{o}(a,m+1,n,k:a_1,\ldots,a_k))}\\
    &+\frac{M_{q}(H_{e}(a,m+1,n+1,k:a_1,\ldots,a_k))}{M_{q}(H_{e}(a,m+1,n,k:a_1,\ldots,a_k))}\cdot\frac{M_{q}(H_{o}(a,m+1,n-1,k:a_1,\ldots,a_k))}{M_{q}(H_{o}(a,m+1,n,k:a_1,\ldots,a_k))}\\
    &=1.
\end{aligned}
\end{equation}
Again, using (2.2), (2.3) and (2.1),
\begin{equation}
\begin{aligned}
    &\frac{M_{q}(H_{e}(a,m+2,n,k:a_1,\ldots,a_k))}{M_{q}(H_{e}(a,m+1,n,k:a_1,\ldots,a_k))}\\
    &
    \begin{aligned}
        =&\prod_{i=1}^{n+k}\Bigg[\frac{q^{(m+k+2)-i}+q^{i-(m+k+2)}}{2}\Bigg]\cdot\frac{\langle m+k+1\rangle_{q}!\langle 2a+2b_1+m+n+2k+1\rangle_{q}!}{\langle m+n+2k+1\rangle_{q}!\langle 2a+2b_1+m+k+1\rangle_{q}!}\\
        &\cdot\prod_{i=0}^{k-1}\Bigg[\frac{\langle a+b_1+b_{i+1}+m+k-i+1\rangle_{q^2}\langle \frac{1}{2}m+\frac{1}{2}n+k-i+\frac{1}{2}\rangle_{q^2}\langle \frac{1}{2}m+\frac{1}{2}n+k-i\rangle_{q^2}}{\langle m+k-i+1\rangle_{q^2}\langle a+b_1+\frac{1}{2}m+\frac{1}{2}n+k-i+\frac{1}{2}\rangle_{q^2}\langle b_{i+1}+\frac{1}{2}m+\frac{1}{2}n+k-i\rangle_{q^2}}\Bigg],    
    \end{aligned}
\end{aligned}
\end{equation}
\begin{equation}
\begin{aligned}
    &\frac{M_{q}(H_{o}(a,m,n,k:a_1,\ldots,a_k))}{M_{q}(H_{o}(a,m+1,n,k:a_1,\ldots,a_k))}\\
    &
    \begin{aligned}
        =&\prod_{i=1}^{n+k}\Bigg[\frac{2}{q^{(m+k+1)-i}+q^{i-(m+k+1)}}\Bigg]\cdot\frac{\langle m+n+2k\rangle_{q}!\langle 2a+2b_1+m+k+1\rangle_{q}!}{\langle m+k\rangle_{q}!\langle 2a+2b_1+m+n+2k+1\rangle_{q}!}\\
        &\cdot\prod_{i=0}^{k-1}\Bigg[\frac{\langle m+k-i\rangle_{q^2}\langle a+b_1+\frac{1}{2}m+\frac{1}{2}n+k-i+\frac{1}{2}\rangle_{q^2}\langle b_{i+1}+\frac{1}{2}m+\frac{1}{2}n+k-i\rangle_{q^2}}{\langle a+b_1+b_{i+1}+m+k-i+1\rangle_{q^2}\langle \frac{1}{2}m+\frac{1}{2}n+k-i-\frac{1}{2}\rangle_{q^2}\langle \frac{1}{2}m+\frac{1}{2}n+k-i\rangle_{q^2}}\Bigg],    
    \end{aligned}
\end{aligned}
\end{equation}
\begin{equation}
\begin{aligned}
    &\frac{M_{q}(H_{e}(a,m+1,n+1,k:a_1,\ldots,a_k))}{M_{q}(H_{e}(a,m+1,n,k:a_1,\ldots,a_k))}\\
    &
    \begin{aligned}
        =&\prod_{i=1}^{m+k+1}\Bigg[\frac{q^{(n+k+1)-i}+q^{i-(n+k+1)}}{2}\Bigg]\cdot\frac{\langle n+k\rangle_{q}!\langle 2a+2b_1+m+n+2k+1\rangle_{q}!}{\langle m+n+2k+1\rangle_{q}!\langle 2a+2b_1+n+k\rangle_{q}!}\\
        &\cdot\prod_{i=0}^{k-1}\Bigg[\frac{\langle a+b_1+b_{i+1}+n+k-i\rangle_{q^2}\langle \frac{1}{2}m+\frac{1}{2}n+k-i+\frac{1}{2}\rangle_{q^2}\langle \frac{1}{2}m+\frac{1}{2}n+k-i\rangle_{q^2}}{\langle n+k-i\rangle_{q^2}\langle a+b_1+\frac{1}{2}m+\frac{1}{2}n+k-i+\frac{1}{2}\rangle_{q^2}\langle b_{i+1}+\frac{1}{2}m+\frac{1}{2}n+k-i\rangle_{q^2}}\Bigg],    
    \end{aligned}
\end{aligned}
\end{equation}
and
\begin{equation}
\begin{aligned}
    &\frac{M_{q}(H_{o}(a,m+1,n-1,k:a_1,\ldots,a_k))}{M_{q}(H_{o}(a,m+1,n,k:a_1,\ldots,a_k))}\\
    &
    \begin{aligned}
        =&\prod_{i=1}^{m+k+1}\Bigg[\frac{2}{q^{(n+k)-i}+q^{i-(n+k)}}\Bigg]\cdot\frac{\langle m+n+2k\rangle_{q}!\langle 2a+2b_1+n+k\rangle_{q}!}{\langle n+k-1\rangle_{q}!\langle 2a+2b_1+m+n+2k+1\rangle_{q}!}\\
        &\cdot\prod_{i=0}^{k-1}\Bigg[\frac{\langle n+k-i-1\rangle_{q^2}\langle a+b_1+\frac{1}{2}m+\frac{1}{2}n+k-i+\frac{1}{2}\rangle_{q^2}\langle b_{i+1}+\frac{1}{2}m+\frac{1}{2}n+k-i\rangle_{q^2}}{\langle a+b_1+b_{i+1}+n+k-i\rangle_{q^2}\langle \frac{1}{2}m+\frac{1}{2}n+k-i-\frac{1}{2}\rangle_{q^2}\langle \frac{1}{2}m+\frac{1}{2}n+k-i\rangle_{q^2}}\Bigg].
    \end{aligned}
\end{aligned}
\end{equation}
Hence,
\begin{equation}
\begin{aligned}
    &\frac{M_{q}(H_{e}(a,m+2,n,k:a_1,\ldots,a_k))}{M_{q}(H_{e}(a,m+1,n,k:a_1,\ldots,a_k))}\cdot\frac{M_{q}(H_{o}(a,m,n,k:a_1,\ldots,a_k))}{M_{q}(H_{o}(a,m+1,n,k:a_1,\ldots,a_k))}\\
    &+\frac{M_{q}(H_{e}(a,m+1,n+1,k:a_1,\ldots,a_k))}{M_{q}(H_{e}(a,m+1,n,k:a_1,\ldots,a_k))}\frac{M_{q}(H_{o}(a,m+1,n-1,k:a_1,\ldots,a_k))}{M_{q}(H_{o}(a,m+1,n,k:a_1,\ldots,a_k))}\\
    &
    \begin{aligned}
        =&\frac{q^{m+k+1}+q^{-(m+k+1)}}{q^{n-m-1}+q^{-(n-m-1)}}\cdot\frac{\langle m+k+1\rangle_{q}}{\langle m+n+2k+1\rangle_{q}}\cdot\frac{\langle m+1\rangle_{q^2}\langle \frac{1}{2}m+\frac{1}{2}n+k+\frac{1}{2}\rangle_{q^2}}{\langle m+k+1\rangle_{q^2}\langle \frac{1}{2}m+\frac{1}{2}n+\frac{1}{2}\rangle_{q^2}}\\
        &+\frac{q^{n+k}+q^{-(n+k)}}{q^{n-m-1}+q^{-(n-m-1)}}\cdot\frac{\langle n+k\rangle_{q}}{\langle m+n+2k+1\rangle_{q}}\cdot\frac{\langle n\rangle_{q^2}\langle \frac{1}{2}m+\frac{1}{2}n+k+\frac{1}{2}\rangle_{q^2}}{\langle n+k\rangle_{q^2}\langle \frac{1}{2}m+\frac{1}{2}n+\frac{1}{2}\rangle_{q^2}}
    \end{aligned}
    \\
    &=\frac{\langle m+1\rangle_{q^2}+\langle n\rangle_{q^2}}{(q^{n-m-1}+q^{-(n-m-1)})\langle \frac{1}{2}m+\frac{1}{2}n+\frac{1}{2}\rangle_{q^2}}\\
    &=1.
\end{aligned}
\end{equation}
and this verifies Step 3$'$.

We just showed that $M_{q}(H_{e}(a,m,n,k:a_1,\ldots,a_{k-1},a_k))$ and $M_{q}(H_{o}(a,m,n,k:a_1,\ldots,a_{k-1},a_k))$ are given by (2.2) and (2.3) when $n-m=E$ and the depths of the regions are $D$. Thus, by inner induction, we can conclude that the TGFs of the two families of regions are always given by (2.2) and (2.3) when $n-m=E$. This verifies the induction step of the outer induction, so we can conclude that (2.2) and (2.3) are true for any $m$ and $n$ with $n\geq m$. This completes the proof of Theorem 2.2. \qedsymbol

\end{document}